\documentclass[12pt]{article}
\usepackage[utf8]{inputenc}
\usepackage{mathtools}
\usepackage{amsmath, amssymb, amsfonts, amsthm}
\usepackage{mathrsfs}
\usepackage[german, american]{babel}

\newtheorem{theorem}{Theorem}
\newtheorem{lemma}[theorem]{Lemma}
\newtheorem{remark}[theorem]{Remark}
\newtheorem{corollary}[theorem]{Corollary}
\newtheorem{definition}[theorem]{Definition}

\allowdisplaybreaks

\newcommand{\opA}{\mathcal{A}}

\newcommand{\opK}{\mathcal{K}}
\newcommand{\opL}{\mathcal{L}}

\newcommand{\spH}{\mathcal{H}}
\newcommand{\p}{\partial_t}
\newcommand{\ii}{\text{i}}

\newcommand{\nlE}{\mathcal{E}^{\mathrm{nlin}}}
\newcommand{\normE}{\bar{\mathcal{E}}^{\mathrm{lin}}}
\newcommand{\lE}{\mathcal{E}^{\mathrm{lin}}}
\newcommand{\mE}{\mathcal{E}^{\mathrm{diff}}}

\newcommand{\oA}{\mathcal{A}}
\newcommand{\oB}{\mathcal{B}}

\allowdisplaybreaks

\begin{document}
\title{Multidimensional Thermoelasticity \\ for Nonsimple Materials -- \\
Well-Posedness and Long-Time Behavior}

\author{
Andrii Anikushyn\thanks{Faculty of Cybernetics, Taras Shevchenko National University of Kyiv, Kyiv, Ukraine \hfill \texttt{anik\_andrii@univ.kiev.ua}} \and
Michael Pokojovy\thanks{Department of Mathematics, Karlsruhe Institute of Technology, Karlsruhe, Germany \hfill \texttt{michael.pokojovy@kit.edu}}
}

\date{\today}

\pagestyle{myheadings}
\thispagestyle{plain}
\markboth{\textsc{A. Anikushyn and M. Pokojovy}}{\textsc{Multidimensional Thermoelasticity for Nonsimple Materials }}

\maketitle

\begin{abstract}
	An initial-boundary value problem for the multidimensional type III thermoelaticity 
	for a nonsimple material with a center of symmetry is considered.
	In the linear case, the well-posedness with and without
	Kelvin-Voigt and/or frictional damping in the elastic part
	as well as the lack of exponential stability in the elastically undamped case is proved.
	Further, a frictional damping for the elastic component is shown to lead to the exponential stability.
	A Cattaneo-type hyperbolic relaxation for the thermal part is introduced
	and the well-posedness and uniform stability under a nonlinear frictional damping are obtained
	using a compactness-uniqueness-type argument.
	Additionally, a connection between the exponential stability
	and exact observability for unitary $C_{0}$-groups is established.
\end{abstract}


{\bf AMS:}
	35B40;   	% Asymptotic behavior of solutions
	35G46;   	% Initial-boundary value problems for linear higher-order systems
	35G61;   	% Initial-boundary value problems for nonlinear higher-order systems
	35Q74;   	% PDEs in connection with mechanics of deformable solids
	74F05;   	% Coupling of solid mechanics with other effects: Thermal effects
	93D20    	% Asymptotic stability

%%%%%%%%%%%%%%%%%%%%%%%%%%%%%%%%%%%%%%%%%%%%%%%%%%%%%%%%%%%%%%%%%%%%%%%%%%%%%%%%%%%%%%%%%%%%%%%%%%%%%%%%%%%%%%%%%%%%%%%%%%
%%%%%%%%%%%%%%%%%%%%%%%%%%%%%%%%%%%%%%%%%%%%%%%%%%%%%%%%%%%%%%%%%%%%%%%%%%%%%%%%%%%%%%%%%%%%%%%%%%%%%%%%%%%%%%%%%%%%%%%%%%

\section{Introduction}
In modern rational mechanics,
when describing a material body,
the general principles such as field equations or jump conditions
are separated from the so-called constitutive equations or material laws
(cf. \cite[pp. 1--2]{TrNo2004}).
Whereas the former are common among all representatives
of a major class of materials such as solids or fluids,
the latter are meant to uniquely characterize each particular material.
Typically, a material law describes the response of a material to various stimuli applied
and is modelled by an algebraic or an operator equation
(see, e.g., \cite[Chapter C3]{TrNo2004}).
In macroscopic theories such as the classical theory of (thermo)elasticity
or its modern generalizations and unifications such as the one proposed
by Green and Naghdi \cite{GrNa1995.I, GrNa1995.II, GrNa1995.III}, etc.,
the microscopic structure of the material is ignored
and the material law is obtained in the form of a stress-strain relation,
which can be measured experimentally.
Recently, Rajagopal \cite{Ra2003} made a strong case
for using implicit material laws, which are physically more sound
as they let the stress, i.e., the force,
induce the strain and not vice versa.
In anelastic bodies, e.g., the thermoelastic ones,
the entropy and the resulting irreversibility play a very important role.
In the theory of thermoelasticity,
a relation between the heat flux and the temperature gradient is postulated
such that the second theorem of thermodynamics,
usually in form of Clausius-Duhem inequality, is satisfied \cite[Sections 96 and 96 bis]{TrNo2004}.
Alternatively, an entropy balance equation can be used
as proposed by Green and Naghdi \cite{GrNa1995.II, GrNa1995.III}.

When being applied to various materials with microstructure,
the macroscopic theory of (thermo)elasticity
does not provide an adequate mechanical and thermodynamical description.
They include but are not limitted to
porous elastic media, micropolar elastic solids, materials with microstructure
and nonsimple elastic solids,
which were first introduced in the works of
Truesdell and Toupin \cite{TruTou1960}, Green and Rivlin \cite{GrRi1964}, Mindlin \cite{Mi1964} and Toupin \cite{Tou1964}.
Their common feature is that they incorprorate extra field variables
such as microstresses or microrotations, hyperstresses or volume fractions, etc.
For further details, we refer the reader
to comprehensive monographs by Ciarletta and Ie\c{s}an \cite{CiaIes1993}
and Ie\c{s}an \cite{Ies2004}.

We now briefly summarize the linear thermoelasticity model
for nonsimple materials without energy dissipation introduced by Quintanilla \cite{Quin1}.
To this end, consider a rigid body occupying in a reference configuration a bounded domain $\Omega$ of $\mathbb{R}^{d}$
with the Lipschitz boundary $\Gamma := \partial \Omega$.
For a period of time $t \geq 0$ and a material point $\mathbf{x} = (x_{i}) \in \Omega$,
let $\mathbf{u} = (u_{i})$, $(u_{iL})$ and $T$, being functions $\mathbf{x}$ and $t$,
denote the displacement vector,
the homogeneous deformation tensor of the `particle' with its center of mass located at $\mathbf{x}$
and the relative temperature measured with respect to a constant reference temperature $T_{0}$,
which is assumed to be attained at some time $t_{0} \geq 0$.
Further, we define the thermal displacement
\begin{equation}
	\tau(\mathbf{x}, t) := \int_{t_{0}}^{t} T(\mathbf{x}, s) \mathrm{d}s. \notag
\end{equation}
With $\mathbf{t} = (t_{Kj})$ denoting the first Piola \& Kirchhoff stress tensor,
the linear balance of momentum reads as
\begin{equation}
	\rho \ddot{u}_{i} = t_{Ki, K} + \rho f_{i}, \notag
\end{equation}
where $\rho$ stands for the material density and $\mathbf{f} = (f_{i})$ is the volumetric force.
Here and the sequel,
we employ the Einstein's summation convention
as well as the standard notation for temporal and spatial
derivatives of scalar and tensor fields (cf. \cite[Chapter 1]{Ies2004}).
In contrast to the classical linear theory of thermoelasticity,
which utilizes the (generalized) Hooke's law
to postulate a linear relation between the elastic part of $\mathbf{t}$
and the infinitesimal Cauchy strain tensor
\begin{equation}
	\boldsymbol{\varepsilon} = \frac{1}{2} (\nabla \mathbf{u}^{T} + (\nabla \mathbf{u})\big), \notag
\end{equation}
the hyperstresses need to be accounted for.
The latter can be shown to incorporate higher order derivatives of $\mathbf{u}$
(cf. \cite[Chapter 7.1]{Ies2004}).
Following Quintanilla \cite[Sections 1 and 2]{Quin1},
the linear equations of thermoelasticity without energy dissipation for a simple material with a center of symmetry read as
\begin{align}
	\rho \ddot{u}_{i} &= \big(A_{iJRs} u_{s, R} - \beta_{Ji} T - (C_{iJKSRl} u_{l, RS} + M_{iJKR} \tau_{, R})_{, K}\big)_{, J} + \rho f_{i}, 
	\label{EQUATION_LINEAR_NONSIMPLE_THERMOELASTICITY_GENERAL_PDE_1} \\
	a \ddot{\tau} &=
	-\beta_{Ki} \dot{u}_{i, K} - M_{jLKI} u_{j, LKI} + K_{IJ} \tau_{, IJ} + \rho T_{0}^{-1} R
	\label{EQUATION_LINEAR_NONSIMPLE_THERMOELASTICITY_GENERAL_PDE_2}
\end{align}
in $\Omega \times (0, \infty)$, where $R$ stands for the volumetric heat sources.
We refer the reader to \cite[Section 2]{Quin1} for an explanation of material constants
in Equations (\ref{EQUATION_LINEAR_NONSIMPLE_THERMOELASTICITY_GENERAL_PDE_1})--(\ref{EQUATION_LINEAR_NONSIMPLE_THERMOELASTICITY_GENERAL_PDE_2})
as well as a discussion on their symmetry and positive definiteness properties.
The homogeneous Dirichlet-Dirchlet boundary conditions for
Equations (\ref{EQUATION_LINEAR_NONSIMPLE_THERMOELASTICITY_GENERAL_PDE_1})--(\ref{EQUATION_LINEAR_NONSIMPLE_THERMOELASTICITY_GENERAL_PDE_2}) 
on $\Gamma$ are given by
\begin{equation}
	u_{i} = 0, \quad u_{i, A} = 0, \quad \tau = 0 \text{ in } \Gamma \times (0, \infty),
	\label{EQUATION_LINEAR_NONSIMPLE_THERMOELASTICITY_GENERAL_BC}
\end{equation}
whereas the initial conditions are stated as
\begin{equation}
	u_{i}(\cdot, 0) = u_{i}^{0}, \quad
	\dot{u}_{i}(\cdot, 0) = u_{i}^{1}, \quad
	\tau(\cdot, 0) = \tau^{0}, \quad
	T(\cdot, 0) = T^{0} \text{ in } \Omega.
	\label{EQUATION_LINEAR_NONSIMPLE_THERMOELASTICITY_GENERAL_IC}
\end{equation}

In \cite{Quin1},
Quintanilla proposed a logarithmically convex energy-like function
to prove the uniqueness for Equations 
(\ref{EQUATION_LINEAR_NONSIMPLE_THERMOELASTICITY_GENERAL_PDE_1})--(\ref{EQUATION_LINEAR_NONSIMPLE_THERMOELASTICITY_GENERAL_IC}).
He further employed the operator semigroup theory to obtain the existence of solutions.

Fern\'{a}ndez Sare et al. \cite{FeSaMuRiQui2010}
considered a one-dimensional counterpart of
Equations (\ref{EQUATION_LINEAR_NONSIMPLE_THERMOELASTICITY_GENERAL_PDE_1})--(\ref{EQUATION_LINEAR_NONSIMPLE_THERMOELASTICITY_GENERAL_IC})
formulated in terms of Green \& Naghdi's type I thermoelasticity (cf. \cite{GrNa1995.II})
under various sets of boundary conditions.
Based on the operator semigroup theory and Gearhart \& Pr\"uss' theorem,
the well-posedness and exponential stability of solution were shown.
Further, the spectral analyticity criterion was employed to show the lack of analyticity for the underlying semigroup.
Finally, the authors proved the impossibility
of the solutions to localize in time.
Similar results have later been obtained
by Maga\~{n}a and Quintanilla \cite{MaQui2014} also for the case of the more comprehensive type III thermoelasticity.

Pata and Quintanilla \cite{PaQui2010}
considered a 3D version of Equation (\ref{EQUATION_LINEAR_NONSIMPLE_THERMOELASTICITY_GENERAL_PDE_1}) with $T \equiv 0$
for the case vanishing anti-plane share deformations (i.e., $u_{1} = u_{2} = 0$).
The Equation for $u := u_{3}$ reduces then to
\begin{equation}
	\rho \ddot{u} = \mu(0) \triangle u
	- l_{2}(0) \triangle^{2} u
	+ \int_{0}^{\infty} \big(\mu'(s) \triangle u(t - s) - l_{2}'(s) \triangle^{2} u(t - s)\big) \mathrm{d}s = 0
	\label{EQUATION_PA_QUI_EQ1}
\end{equation}
together with the boundary condition
\begin{equation}
	u = \triangle u = 0 \label{EQUATION_PA_QUI_EQ2}
\end{equation}
and the initial condition
\begin{equation}
	u(\cdot, -s) = g(\cdot, s) \text{ for } s \in [0, \infty). \label{EQUATION_PA_QUI_EQ3}
\end{equation}
First, an abstract version of the initial-boundary value problem
(\ref{EQUATION_PA_QUI_EQ1})--(\ref{EQUATION_PA_QUI_EQ3}) was studied using the operator semigroup theory.
Based on a modification of Gearhart \& Pr\"uss' theorem,
a condition for the exponential stability was derived.

Gawinecki and \L{}azuka \cite{GaLu2006} considered a Cauchy problem
for the genuinely nonlinear version of
Equations (\ref{EQUATION_LINEAR_NONSIMPLE_THERMOELASTICITY_GENERAL_PDE_1})--(\ref{EQUATION_LINEAR_NONSIMPLE_THERMOELASTICITY_GENERAL_PDE_2})
within Green \& Naghdi's type I thermoelasticity for homogeneous isotropic media.
Under appropriate conditions on the nonlinearity,
a global classical solutions was obtained based on $L^{p}$-$L^{q}$-estimates.

Aouadi and Moulahi \cite{AouMou2015} studied an initial-boundary value problem
for the anti-plane share component of a nonsimple thermoelastic body
with a control distributed over an open subset $\omega$ of $\Omega$
\begin{align}
	u_{tt} - c^{2} \triangle u + \alpha \triangle^{2} u + c^{2} \gamma (-\triangle)^{1/2} \theta &= \chi_{\omega} u_{1}, 
	\label{EQUATION_AOUADI_MOULAHI_PDE_1} \\
	\theta_{t} - \triangle \theta - \gamma (-\triangle)^{1/2} u_{t} &= \chi_{\omega} u_{2}
	\label{EQUATION_AOUADI_MOULAHI_PDE_2}
\end{align}
subject to the homogeneous Dirichlet-Dirichlet boundary conditions
and the standard initial conditions. Here, $\chi_{\omega}$ stands for the indicator function of $\omega$.
For the uncontrolled case, i.e., $u_{1} = u_{2} \equiv 0$,
a well-posedness and exponential stability result for
Equations (\ref{EQUATION_AOUADI_MOULAHI_PDE_1})--(\ref{EQUATION_AOUADI_MOULAHI_PDE_2}) was proved.
Further, the authors showed
the system is approximately controllable
over the standard state space
by $L^{2}\big(0, T; L^{2}(\Omega)\big)$-controls at any time $T > 0$
if $\omega$ is nonempty.

In the present paper, we first consider Equations
(\ref{EQUATION_LINEAR_NONSIMPLE_THERMOELASTICITY_GENERAL_PDE_1})--(\ref{EQUATION_LINEAR_NONSIMPLE_THERMOELASTICITY_GENERAL_IC})
restated within the type III thermoelasticity
(cf. \cite[Section 2]{MaQui2014} for the 1D system)
with a (macroscopic) Kelvin-Voigt damping for $\mathbf{u}$.
Hence, our equations read as
% \begin{align}
% 	\label{SystemWithoutFriction1}
% 	\rho \ddot{u}_i =& \left( A_{iJKj} u_{j,K} - \beta_{Ji} \theta -
% 	(C_{iJKLIj} u_{j,IL} + M_{iJKL} \tau_{,L})_{,K} \right)_{,J} \\  
% 	%
% 	\label{SystemWithoutFriction2}
% 	a \ddot{\tau} =&
% 	-\beta_{Ki} \dot{u}_{i,K} + M_{jLKI} u_{j,LKI} + K_{IJ}\tau_{,IJ}.
% \end{align}
\begin{align}
	\label{EQUATION_NONSIMPLE_THERMOELASTICITY_LINEAR_KELVIN_VOIGT_DAMPING_1}
	\rho \ddot{u}_i &= \big(A_{iJKj} u_{j,K} - \beta_{Ji} \dot{\tau} -
	(C_{iJKLIj} u_{j,IL} + M_{iJKL} \tau_{,L} - B_{iJKj}\dot{u}_{j})_{,K}\big)_{,J} \\
	\label{EQUATION_NONSIMPLE_THERMOELASTICITY_LINEAR_KELVIN_VOIGT_DAMPING_2}
	a \ddot{\tau} &=
	-\beta_{Ki} \dot{u}_{i,K} + m_{IJ}\dot{\tau}_{,IJ} + M_{jLKI} u_{j,LKI} + K_{IJ}\tau_{,IJ}.
\end{align}
together with the boundary and initial conditions in Equations
(\ref{EQUATION_LINEAR_NONSIMPLE_THERMOELASTICITY_GENERAL_BC})--(\ref{EQUATION_LINEAR_NONSIMPLE_THERMOELASTICITY_GENERAL_IC}).
Throughout the paper,
the following (natural) positivity and positive definitess conditions are postulated:
\begin{itemize}
	\item[I.] $a,\rho > 0$.
	
	\item[II.] There exists $\alpha > 0$ such that
	\begin{align}
		\notag
		C_{iIJKLj} u_{i,JI} u_{j,LK} + M_{iJKL} u_{i,KJ} \tau_{, L} & \\ 
		\label{DPositDefined1}
		+ M_{jLKI} u_{j,LK}\tau_{,I} + K_{IJ} \tau_{,I}\tau_{,J} & \ge
		\alpha(\tau_{,R}\tau_{,R} + u_{k,ST}u_{k,ST}), \\ 
		\label{DPositDefined2} A_{iKLj} u_{i,K} u_{j,L} & \ge \alpha
		u_{k,R}u_{k,R}.
	\end{align}
\end{itemize}
In Section \ref{SECTION_WELL_POSEDNESS}, we show
Equations (\ref{EQUATION_NONSIMPLE_THERMOELASTICITY_LINEAR_KELVIN_VOIGT_DAMPING_1})--(\ref{EQUATION_NONSIMPLE_THERMOELASTICITY_LINEAR_KELVIN_VOIGT_DAMPING_2}),
(\ref{EQUATION_LINEAR_NONSIMPLE_THERMOELASTICITY_GENERAL_BC})--(\ref{EQUATION_LINEAR_NONSIMPLE_THERMOELASTICITY_GENERAL_IC}) are well-posed.
Note that no positive definiteness conditions on $(B_{iJKj})$ or $(m_{IJ})$ are imposed.
By a standard perturbation argument (e.g., \cite[Chapter 3.1, Theorem 1.1]{Pazy1},
the problem remains well-posed if the Kelvin-Voigt damping
is replaced or complemented with a frictional damping.

In Section \ref{SECTION_LACK_OF_EXPONENTIAL_STABILITY}, we prove that, 
in contrast to the 1D situation (cf. \cite[Section 2]{MaQui2014}),
if $B_{iJKj} \equiv 0$,
the positive definiteness of $(m_{IJ})$ is not sufficient
to exponentially stabilize the system.
Motivated by the necessity of an additional damping for $\mathbf{u}$,
we let $B_{iJKj} \equiv 0$, but add a linear frictional damping
to Equation (\ref{EQUATION_NONSIMPLE_THERMOELASTICITY_LINEAR_KELVIN_VOIGT_DAMPING_1}) thus obtaining
\begin{equation}
	\label{EQUATION_NONSIMPLE_THERMOELASTICITY_LINEAR_FRICTIONAL_DAMPING}
	\rho \ddot{u}_i = \big(A_{iJKj} u_{j,K} - \beta_{Ji} \dot{\tau} -
	(C_{iJKLIj} u_{j,IL} + M_{iJKL} \tau_{,L} - B_{iJKj}\dot{u}_{j})_{,K}\big)_{,J} - E_{ij} \dot{u}_{j}.
\end{equation}
Further, in Section \ref{SECTION_EXPONENTIAL_STABILITY},
under a positive definiteness assumption on $(E_{ij})$,
we use the Lyapunov's method to show the exponential stability of Equations
(\ref{EQUATION_NONSIMPLE_THERMOELASTICITY_LINEAR_FRICTIONAL_DAMPING}), (\ref{EQUATION_NONSIMPLE_THERMOELASTICITY_LINEAR_KELVIN_VOIGT_DAMPING_2}),
(\ref{EQUATION_LINEAR_NONSIMPLE_THERMOELASTICITY_GENERAL_BC})--(\ref{EQUATION_LINEAR_NONSIMPLE_THERMOELASTICITY_GENERAL_IC}).

In Section \ref{SECTION_HYPERBOLIZED_SYSTEM}, a hyperbolization of (\ref{EQUATION_NONSIMPLE_THERMOELASTICITY_LINEAR_KELVIN_VOIGT_DAMPING_2})
together a semilinear frictional damping in Equation (\ref{EQUATION_NONSIMPLE_THERMOELASTICITY_LINEAR_KELVIN_VOIGT_DAMPING_1}) are considered.
The resulting semilinear system reads then as
\begin{align}
	\rho \ddot{u}_i &= \left( A_{iJKj} u_{j,K} - \beta_{Ji} \dot{\tau} -
	(C_{iJKLIj} u_{j,IL} + M_{iJKL} \tau_{,L} )_{,K} \right)_{,J} -
	E(|\dot{u}|) \dot{u}_i,
	\label{EQUATION_NONSIMPLE_THERMOELASTICITY_SEMILINEAR_DAMPED_HYPERBOLIZED_1} \\  
	a \ddot{\tau} &=
	-\beta_{Ki} \dot{u}_{i,K} + m_{IJ}q_{I,J} + M_{jLKI} u_{j,LKI} +
	K_{IJ}\tau_{,IJ},
	\label{EQUATION_NONSIMPLE_THERMOELASTICITY_SEMILINEAR_DAMPED_HYPERBOLIZED_2} \\
	\kappa \dot{q}_i &= \dot{\tau}_{,i} - q_{i},
	\label{EQUATION_NONSIMPLE_THERMOELASTICITY_SEMILINEAR_DAMPED_HYPERBOLIZED_3}
\end{align}
with a small parameter $\kappa > 0$.
Here, $(m_{IJ} q_{I, J})$ plays the role of the heat flux.
If the function $E$ is constant,
similar to the case $\kappa = 0$,
Equations (\ref{EQUATION_NONSIMPLE_THERMOELASTICITY_SEMILINEAR_DAMPED_HYPERBOLIZED_1})--(\ref{EQUATION_NONSIMPLE_THERMOELASTICITY_SEMILINEAR_DAMPED_HYPERBOLIZED_3})
are exponentially stable.
In the nonlinear case,
we exploit a technique due to Lasiecka and Tataru \cite{Lasiecka1}
to prove the uniform stability of the nonlinear system
dependent on the behavior of $E$ at $0$ and infinity.
Here, we use a generalization of a technique dating back to Haraux \cite{Haraux1},
which is given in the appendix \ref{SECTION_APPENDIX}.

%%%%%%%%%%%%%%%%%%%%%%%%%%%%%%%%%%%%%%%%%%%%%%%%%%%%%%%%%%%%%%%%%%%%%%%%%%%%%%%%%%%%%%%%%%%%%%%%%%%%%%%%%%%%%%%%%%%%%%%%%%
%%%%%%%%%%%%%%%%%%%%%%%%%%%%%%%%%%%%%%%%%%%%%%%%%%%%%%%%%%%%%%%%%%%%%%%%%%%%%%%%%%%%%%%%%%%%%%%%%%%%%%%%%%%%%%%%%%%%%%%%%%

\section{Linear sytem with a Kelvin-Voigt damping for the elastic part: Well-posedness}
\label{SECTION_WELL_POSEDNESS}
In this section, we consider the equations of type III thermoelasticity for a nonsimple material
with a Kelvin-Voigt damping for $\mathbf{u}$.
By a standard perturbation argument \cite[Chapter 3.1, Theorem 1.1]{Pazy1},
the results remain valid also when a frictional damping, i.e., a term like $E_{ij} \dot{u}_{j}$,
is considered instead of or in addition to $B_{iJKj} \dot{u}_{j, KJ}$.
The equations read then as
\begin{align} 
	\label{SystemTwoFrictionSecondCase1}
	\rho \ddot{u}_i &= \big(A_{iJKj} u_{j,K} - \beta_{Ji} \dot{\tau} -
	(C_{iJKLIj} u_{j,IL} + M_{iJKL} \tau_{,L} - B_{iJKj}\dot{u}_{j})_{,K}\big)_{,J} \\
	\label{SystemTwoFrictionSecondCase2}
	a \ddot{\tau} &=
	-\beta_{Ki} \dot{u}_{i,K} + m_{IJ}\dot{\tau}_{,IJ} + M_{jLKI} u_{j,LKI} + K_{IJ}\tau_{,IJ}
\end{align}
in $\Omega \times (0, \infty)$ together with the boundary conditions
\begin{equation}
	\label{WPBoundaryConditions}
	u_i = 0, \quad u_{i,J} = 0, \quad \tau = 0 \text{ in } \Gamma \times (0, \infty)
\end{equation}
and the initial conditions
\begin{equation}
	\label{WPInitialConditions}
	u_i(\cdot, 0) = u_i^0, \quad \dot{u}_i(\cdot,0) = \dot{u}_i^0, \quad 
	\tau (\cdot, 0) = \tau^0, \quad \dot{\tau}(\cdot, 0) = \dot{\tau}^0 \text{ in } \Omega,
\end{equation}
where $\Gamma := \partial \Omega$ is assumed Lipschitzian.
In addition to conditions I, II, we assume
\begin{align}
	B_{iJKj} \xi_{j,K} \xi_{i,J} \ge 0, \quad
	m_{IJ} \xi_{I} \xi_{J} \ge 0.
	\label{EQUATION_POSITIVE_SEMIDEFINITENESS_DAMPING}
\end{align}

Letting $v = \dot{u}$ and $\theta = \dot{\tau}$, we rewrite the system
(\ref{SystemTwoFrictionSecondCase1})--(\ref{SystemTwoFrictionSecondCase2}) as follows:
\begin{align}
	\label{WPSystem1}
	\dot{u}_i &= v_i, \\ 
	\label{WPSystem2}
	\dot{v}_i &= \frac{1}{\rho} \big(A_{iJKj} u_{j,K} - \beta_{Ji}
	\theta - (C_{iJKLIj} u_{j,IL} + M_{iJKL} \tau_{,L} - B_{iJKj}v_{j})_{,K} \big)_{,J}, \\ 
	\label{WPSystem3}
	\dot{\tau} &= \theta, \\ 
	\label{WPSystem4}
	\dot{\theta} &= \frac{1}{a}  \big(-\beta_{Ki} v_{i,K} +
	(m_{IJ}\theta_{,J} + M_{jLKI} u_{j,LK} + K_{IJ}\tau_{,J})_{,I}\big).
\end{align}
We assume the evolution is taking place on the Hilbert space
\begin{equation}
	\mathcal{H} = \left\{U \,|\, U = (u, v, \tau, \theta)^{T} \in 
	\left(H^2_{0}(\Omega)\right)^{d} \times \left( L^2(\Omega) \right)^{d} \times
	H_0^1(\Omega) \times L^2(\Omega)\right\}, \notag
\end{equation}
where
\begin{equation}
	H^{s}_{0}(\Omega) :=
	\mathrm{clos}\big(C_{0}^{\infty}(\Omega), \|\cdot\|_{H^{s}}\big) \text{ for } s \in \mathbb{N}, \notag
\end{equation}
equipped with the scalar product
\begin{align}
	\langle U, U^{\ast} \rangle_{\mathcal{H}} &= \int_{\Omega} \big(\rho v_{i} v_{i}^* +
	a \theta \theta^* + A_{iKLj} u_{i,K} u^*_{j,L} + C_{iIJKLj} u_{i,JI}
	u^*_{j,LK} + K_{IJ} \tau_{,J} \tau^*_{,I} \big) \mathrm{d}x \notag \\
	&+ \int_{\Omega} \big(M_{iJKL} u_{i,KJ} \tau^*_{,L} + M_{jLKI}
	u^*_{j,KL} \tau_{,I} \big) \mathrm{d}x, \label{EQUATION_SCALAR_PRODUCT_WEIGHTED}
\end{align}
for $U = (u, v, \tau, \theta)^{T}$, $U^{\ast} = (u^{\ast}, v^{\ast}, \tau^{\ast}, \theta^{\ast})^{T} \in \spH$.
Here and in the following, for the sake of simplicity, 
we drop the physical convention to bold the vector and tensor fields.
It is easy to verify that this scalar product is equivalent with the usual product 
induced by the product topology.

Consider the linear operator $\opA: D(\opA) \subset \mathcal{H} \to
\mathcal{H}$ defined by
\begin{equation}
	\label{WPOpADefine}
	\opA U =
	{\small
	\begin{pmatrix}
		v_i,
		\\
		\frac{1}{\rho} \big( A_{iJKj} u_{j,K} - \beta_{Ji} \theta -
		(C_{iJKLIj} u_{j,IL} + M_{iJKL} \tau_{,L} - B_{iJKj}v_{j} )_{,K}
		\big)_{,J},
		\\
		\theta,
		\\
		\frac{1}{a}  \left( -\beta_{Ki} v_{i,K} + (m_{IJ}\theta_{,J} +
		M_{jLKI} u_{j,LK} + K_{IJ}\tau_{,J})_{,I} \right)
	\end{pmatrix}},
\end{equation}
where
\begin{align}
	\begin{split}
	D(\opA) = \big\{ U \in \mathcal{H} \,|\, \phantom{-}&v_i\in H^2_0(\Omega),\
	\theta \in H^1_0(\Omega), \\
	\phantom{-}&(A_{iJKj} u_{j,K} - \beta_{Ji} \theta)_{,J}  \\
	-&(C_{iJKLIj} u_{j,IL} M_{iJKL} \tau_{,L} - B_{iJKj}v_{j} )_{,K} \big)_{,J} \in L^2(\Omega), \\
	\phantom{-}&(m_{IJ}\theta_{,J} + M_{jLKI} u_{j,LK} +
	K_{IJ}\tau_{,J})_{,I} \in L^2(\Omega)\big\}
	\end{split}
	\label{WPOperatorDomain}
\end{align}
Thus, the abstract form of Equations (\ref{WPSystem1})--(\ref{WPSystem4}) is given by
\begin{equation}
	\label{WPCauchyProblem}
	\dot{U}(t) = \opA U(t) \text{ for } t > 0, \quad U(0) = U_0.
\end{equation}
with $U_{0} := (u^{0}_{i}, u^{1}_{i}, \tau^{0}, \tau^{1})^{T}$.
We first prove the following auxiliary lemma.

\begin{lemma}\label{WPMainLemma}
	The operator $\opA$ defined in Equations
	(\ref{WPOpADefine})--(\ref{WPOperatorDomain}) is the infinitesimal generator of
	a $C_0$-semigroup of contractions on $\mathcal{H}$.
\end{lemma}

\begin{proof}
The proof of Lemma \ref{WPMainLemma} resembles the one of \cite[Theorem 2]{Po2014}.

{\it Denseness: }
Utilizing the fact that
\begin{equation}
	\left(C_{0}^{\infty}(\Omega) \right)^d \times
	\left(C_{0}^{\infty}(\Omega)\right)^d \times C_{0}^{\infty}(\Omega)
	\times C_{0}^{\infty}(\Omega) \notag
\end{equation}
is a subset of the domain of $\opA$, we conclude that $D(\opA)$ is
dense in $\mathcal{H}$. 

{\it Disipativity: }
A straightforward calculation involving the
Green's formula yields
\begin{equation}
	\langle \opA U , U \rangle_{\mathcal{H}} = - \int_{\Omega} \big(B_{iJKj}
	v_{i,K} v_{i,J} + m_{IJ} \theta_{,I} \theta_{,J} \big) \mathrm{d}x.
\end{equation}
Using Equation (\ref{EQUATION_POSITIVE_SEMIDEFINITENESS_DAMPING}), we get
\begin{equation}
	\langle \opA U , U \rangle_{\mathcal{H}} \le -\alpha \int_{\Omega}
	\big(v_{i,K}v_{i,K} + \theta_{,I}\theta_{,I} \big) \mathrm{d}x \le 0.
\end{equation}

{\it Maximality: }
Now we prove that $0\in \rho(\opA)$ with $\rho(\opA)$ standing for
the resolvent set of $\opA$. For $F = (F^1, F^2, F^3, F^4)^{T} \in \mathcal{H}$, 
consider the operator equation $\opA U = F$ or, expicitly,
\begin{align*}
	v_i &= F^1_i,
	\\
	\frac{1}{\rho} \left( A_{iJKj} u_{j,K} - \beta_{Ji} \theta -
	(C_{iJKLIj} u_{j,IL} + M_{iJKL} \tau_{,L} - B_{iJKj}v_{j} )_{,K}
	\right)_{,J} &= F^2_i,
	\\
	\theta &= F^3,
	\\
	\frac{1}{a}  \left( -\beta_{Ki} v_{i,K} + (m_{IJ}\theta_{,J} +
	M_{jLKI} u_{j,LK} + K_{IJ}\tau_{,J})_{,I} \right) &= F^4.
\end{align*}
Eliminating $v_i$ and $\theta$, we obtain
\begin{align}
	\label{WPLaxMilgramProblem1}
	\frac{1}{\rho} \left( A_{iJKj} u_{j,K} - \beta_{Ji} F^3 -
	(C_{iJKLIj} u_{j,IL} + M_{iJKL} \tau_{,L} - B_{iJKj}F^1_{j} )_{,K}
	\right)_{,J} &= F^2_i,
	\\
	\label{WPLaxMilgramProblem2}  \frac{1}{a}  \left( -\beta_{Ki}
	F^1_{i,K} + (m_{IJ}F^3_{,J} + M_{jLKI} u_{j,LK} +
	K_{IJ}\tau_{,J})_{,I} \right) &= F^4.
\end{align}
To solve this system, we exploit the lemma of Lax \& Milgram.
We consider the Hilbert space
\begin{equation}
	\mathcal{V} = \left(H_0^2(\Omega)\right)^d \times H_0^1(\Omega) \notag
\end{equation}
equipped with the standard inner product associated with the product
topology and introduce the bilinear form $\frak{a}: \mathcal{V}
\times \mathcal{V} \to \mathbb{R}$ via
\begin{align*}
	\frak{a}(V,V^*) &= \int_{\Omega} (A_{iJKj} u_{s,R} u^*_{i,J} +
	C_{iIJKLj}u_{j,LK} u^*_{i,IJ}) \ \mathrm{d}x \\
	&+\int_{\Omega} (M_{iJKL} \tau_{,L} u^*_{i,JK} + M_{jLKI} u_{j,KL}
	\tau^*_{,I} + K_{IJ} \tau_{,I} \tau^*_{,J})\ \mathrm{d}x,
\end{align*}
where $V=(u, \tau)^{T}$ and $V^* = (u^*, \tau^*)^{T}$. 
After multiplying Equations (\ref{WPLaxMilgramProblem1})--(\ref{WPLaxMilgramProblem2})
in the inner product $(L^2(\Omega))^n$ and $L^2(\Omega)$ with $\rho W^1$ and $a W^2$, respectively, 
summing up the resulting equations and integrating by parts, 
we obtain a weak formulation of Equations (\ref{WPLaxMilgramProblem1})--(\ref{WPLaxMilgramProblem2}) in the form: 
\begin{align}
	\begin{split}
		&\text{Determine $V\in\mathcal{V}$ such that } \\
		&\frak{a}(V,W) = -\langle G^1, W^1\rangle_{(L^2(\Omega))^n} 
		- \langle G^2, W^2\rangle_{L^2(\Omega)} 
		- m_{IJ}\langle F^3_{,J}, W^2_{,I}\rangle_{L^2(\Omega)}.
	\end{split}
	\label{WPWeakProblem}
\end{align}
Here, $G^1_i = \rho F^2_i + \beta_{Ji} F^3_{, J} - B_{iJKj} F^1_{j,KJ},\
G^2 = aF^4 + \beta_{Ki}F^1_{i,K}$.

The bilinear form $\frak{a}$ is continuous and coercive on $\mathcal{V}$ due to the conditions in
(\ref{DPositDefined1})--(\ref{DPositDefined2}). 
Clearly, the linear functional
\begin{equation}
	(W^1,W^2) \mapsto \langle G^1, W^1\rangle_{(L^2(\Omega))^n} 
	+ \langle G^2, W^2\rangle_{L^2(\Omega)} 
	+ m_{IJ}\langle F^3_{,J}, W^2_{,I}\rangle_{L^2(\Omega)}
	\notag
\end{equation}
is continuous on $\mathcal{V}$. 
Applying now the Lemma of Lax \& Milgram, 
we deduce existence of a unique solution $V = (u, \tau)^{T} \in\mathcal{V}$ 
to Equations (\ref{WPLaxMilgramProblem1})--(\ref{WPLaxMilgramProblem2}). 
Recalling the definition of $v_{i}$ and $\theta$,
we further get $v_i = F^1_i \in H^2_0(\Omega)$ 
and $\theta = F^3 \in H^1_0(\Omega)$. 
Hence, the conditions in Equations (\ref{WPOperatorDomain}) are satisfied. 
Therefore, we have $(u,v,\tau,\theta)^{T} \in D(\opA)$ implying that $(u,v,\tau,\theta)^{T}$ is a strong solution.

Finally we show the continuous dependence of the solution on $F$.
It follows from the Lemma of Lax \& Milgram that there exists $c_0 > 0$ such that
\begin{equation}
	\label{WPAdditionalIneq1}
	\| V \|_{\mathcal{V}} \le c_0 \| \widehat{F} \|_{\mathcal{V}^*},
\end{equation}
where $\widehat{F} = (F^2, F^4)^{T}$. 
Since $\mathcal{V}$ is continuously embedded into $\left(L^2(\Omega)\right)^{d+1}$,  we have
\begin{equation}\label{WPAdditionalIneq2}
	\|\widehat{F}\|_{\mathcal{V}^*} \le c_1 \|\widehat{F}\|_{\left(L^2(\Omega)\right)^{d+1}}.
\end{equation}
Suppose $\|F_{n}\|_{\mathcal{H}} \to 0$ as $n \to \infty$.
By definition,
$\|F^1_{n}\|_{\left(H^2_0(\Omega) \right)^d} \to 0$ and $\|F^3_{n}\|_{H^1_0(\Omega)} \to 0$ as $n \to \infty$ 
and, therefore, $\big\|(v_i)_{n}\big\|_{L^2(\Omega)} \to 0$ and $\|\theta_{n}\|_{L^2(\Omega)} \to 0$ as $n \to \infty$.
Further, $\|F^2_{n}\|_{\left(L^2(\Omega)\right)^d}, \|F^4_{n}\|_{L^2(\Omega)} \to 0$ as $n \to \infty$ 
and it follows from Equations (\ref{WPAdditionalIneq1})--(\ref{WPAdditionalIneq2}) that
$\|V\|_{\mathcal{V}} \to 0$. 
Therefore, the solution $U_{n}$ of equation
$\opA U_{n} = F_{n}$ tends to $0$ in $\mathcal{H}$ as $n \to \infty$.

Now, after we have shown $\opA$ is $m-$dissipative, 
the claim of the present theorem follows from the Lumer \& Phillips theorem 
(cf. \cite[Theorem 4.3, p. 14]{Pazy1}).
\end{proof}

Now, by virtue of \cite[Theorem 1.3, p. 102]{Pazy1},
Equations (\ref{SystemTwoFrictionSecondCase1})--(\ref{WPInitialConditions}) 
as well as their abstract formulation (\ref{WPCauchyProblem}) are well-posed.
\begin{theorem}
	Let $U_0 \in \mathcal{H}$.
	There exists then a unique mild solution $U \in C^0\big([0, \infty), \spH\big)$
	to Equation (\ref{WPCauchyProblem}).
	If $U_0\in D(\opA)$, the mild solution is even a
	classical one satisfying $U \in C^1\big([0, \infty), \spH) \cap C^{0}\big([0, \infty), D(\mathcal{A})\big)$.
\end{theorem}

Applying \cite[Chapter 3.1, Theorem 1.1]{Pazy1}, we further obtain
\begin{corollary}
	Replacing Equation (\ref{SystemTwoFrictionSecondCase1}) with
	\begin{equation}
		\rho \ddot{u}_i = \big(A_{iJKj} u_{j,K} - \beta_{Ji} \dot{\tau} -
		(C_{iJKLIj} u_{j,IL} + M_{iJKL} \tau_{,L} - B_{iJKj}\dot{u}_{j})_{,K}\big)_{,J} +
		E_{ij} \dot{u}_{j}
		\label{SystemTwoFrictionSecondCase1_FRICTIONAL_DAMPING}
	\end{equation}
	for an arbitrary matrix $(E_{ij}) \in \mathbb{R}^{d \times d}$
	in Equations (\ref{SystemTwoFrictionSecondCase1_FRICTIONAL_DAMPING}),
	(\ref{SystemTwoFrictionSecondCase2})--(\ref{WPInitialConditions}),
	the resulting abstract Cauchy problem is well-posed on $\mathcal{H}$.	
\end{corollary}

%%%%%%%%%%%%%%%%%%%%%%%%%%%%%%%%%%%%%%%%%%%%%%%%%%%%%%%%%%%%%%%%%%%%%%%%%%%%%%%%%%%%%%%%%%%%%%%%%%%%%%%%%%%%%%%%%%%%%%%%%%
%%%%%%%%%%%%%%%%%%%%%%%%%%%%%%%%%%%%%%%%%%%%%%%%%%%%%%%%%%%%%%%%%%%%%%%%%%%%%%%%%%%%%%%%%%%%%%%%%%%%%%%%%%%%%%%%%%%%%%%%%%

\section{Linear system with an undamped elastic part: Lack of exponential stability}
\label{SECTION_LACK_OF_EXPONENTIAL_STABILITY}
In this section, we consider Equations
(\ref{SystemTwoFrictionSecondCase1})--(\ref{SystemTwoFrictionSecondCase2}) for the case $B_{iJKj} \equiv 0$, i.e.,
\begin{align} 
	\label{SystemOneFriction1}
	\rho \ddot{u}_i &= \big(A_{iJKj} u_{j,K} - \beta_{Ji} \dot{\tau} -
	(C_{iJKLIj} u_{j,IL} + M_{iJKL} \tau_{,L}, \\
	\label{SystemOneFriction2}
	a \ddot{\tau} &=
	-\beta_{Ki} \dot{u}_{i,K} + m_{IJ}\dot{\tau}_{,IJ} + M_{jLKI} u_{j,LKI} + K_{IJ}\tau_{,IJ}.
\end{align}
Equations (\ref{SystemOneFriction1})--(\ref{SystemOneFriction2})
are precisely the type III thermoelasticity for nonsimple materials.
Under a suitable choice of natural boundary conditions,
i.e., the ones appearing in the Green's formula,
and a feasible selection of the coefficient tensors 
$(A_{iJKj})$, $(\beta_{Ji})$, $(C_{iJKLIj})$, $(m_{IJ})$, $(K_{IJ})$
as well as restricting $\Omega$ to the rectangular configuration
\begin{equation}
	\Omega = (0, \pi) \times (0, \pi), \notag
\end{equation}
we show that the system
(\ref{SystemOneFriction1})--(\ref{SystemOneFriction2}) is lacking exponential stability.
Selecting the relatively open in $\Gamma$ disjunctive sets
\begin{align*}
	\Gamma_{1} := (0, \pi) \times \{0, \pi\} \text{ and }
	\Gamma_{2} := \{0, \pi\} \times (0, \pi),
\end{align*}
we have $\Gamma = \overline{\Gamma}_1 \cup \overline{\Gamma}_2$.
The coefficient tensors
$(A_{iJKj}), (\beta_{Ji}), (C_{iJKLIj}), (m_{IJ}), (K_{IJ})$
are chosen in \cite[p. 5]{Quin1} to model an isotropic material with a center of symmetry.
Further, we select
\begin{equation}
	(m_{ij}) :=
	\left(
	\begin{array}{cc}
		\delta & 0 \\
		0 & \delta \\
	\end{array}\right) \text{ for some } \delta > 0. \notag
\end{equation}
The balance equations and material laws (cf. \cite{Quin1}) read then as
\begin{align*} 
	\tau_{JI,J} &= \sigma_{RJI,RJ} + \rho \ddot{u}_{I},
	\\
	(\phi_{I})_{,I} &= \beta \dot{e}_{LL} - D_{IJ} \dot{\tau}_{,IJ} + a\ddot{\tau},
	\\
	\tau_{IJ} &= \lambda \delta_{IJ} e_{LL} + 2\mu e_{IJ} - \beta \delta_{IJ} \dot{\tau},
	\\
	\sigma_{IJK} &= \frac12 a_1 (\delta_{JK} \kappa_{LLI} + 2\delta_{IJ} \kappa_{KLL} + \delta_{IK} \kappa_{LLJ}) + a_2(
	\delta_{JI} \kappa_{ILL} + \delta_{IK} \kappa_{JLL})
	\\
	&+ 2a_3 \delta_{IJ} \kappa_{LLK} + 2a_4\kappa_{IJK} + a_5(\kappa_{KJI} + \kappa_{KIJ}) + m_1\delta_{IJ} \tau_{,K}
	\\ %\label{NESPhysSystem4}
	&+  m_2(\delta_{JK} \tau_{,I} + \delta_{IK} \tau_{,J}),
	\\
	\phi_{I} &= m_1\kappa_{LLI} + m_2(\kappa_{ILL} + \kappa_{LIL}) +
	\delta_{IJ} \tau_{,J},
	\\
	e_{IJ} &= \frac12 (u_{I,J} + u_{J,I}),\ \kappa_{IJK} = u_{K,IJ}.
\end{align*}
After plugging the constitutive equations into the balance equations, 
we arrive at the system
\begin{align}
	\notag
	\rho \ddot{u}_{1} &= (\lambda+2\mu)u_{1,11} + \mu u_{1,22}
	+ (\lambda + \mu) u_{2,12} - \beta \dot{\tau}_{,1} - b_1 u_{1,1111}
	- b_2 u_{1,1122}
	\\ \label{NESElemSyst1}
	&- b_3 u_{1,2222} - b_4 u_{2,1112} - b_5
	u_{2,1222} -(m_1+2m_2)(\tau_{,111} + \tau_{,122}),
	\\
	\notag
	\rho \ddot{u}_{2} &= (\lambda+2\mu)u_{2,22} + \mu u_{2,11} +
	(\lambda + \mu) u_{1,12} - \beta \dot{\tau}_{,2} - b_1 u_{2,2222} -
	b_2 u_{2,1122}
	\\
	\label{NESElemSyst2} &- b_3 u_{2,1111} - b_4 u_{1,1112} - b_5
	u_{1,1222} -(m_1+2m_2)(\tau_{,112} + \tau_{,222}),
	\\
	\notag
	a\ddot{\tau} &= (m_1 + 2m_2)(u_{1,111}+u_{1,122} + u_{2,112}
	+ u_{2,222}) + \tau_{,1} + \tau_{,2}
	\\
	\label{NESElemSyst3} &- \beta (\dot{u}_{1,1} + \dot{u}_{2,2}) +
	\delta (\dot{\tau}_{,11} + \dot{\tau}_{,22}),
\end{align}
where
\begin{align*}
	b_1 &=  2(a_1+a_2+a_3+a_4+a_5), & b_2 &=  2(a_1+a_2+2a_3+2a_4+a_5), \\
	b_3 &=  2(a_3+a_4),             & b_4 &=  b_5=2(a_1+a_2+a_5).
\end{align*}
Throughout this section, we employ the following set of boundary conditions, 
which naturally arise from partially integrating of the system (\ref{NESElemSyst1})--(\ref{NESElemSyst3}):
\begin{align} \notag
	\sigma_{111} |_{\Gamma}  =
	a_1(u_{1,11}+u_{1,22}+u_{1,12}+u_{2,12}) +2a_2(u_{1,11}+u_{2,21}) &
	\\
	\label{NESBoundaryConditions1}
	+2a_3(u_{1,11} + u_{1,22}) + 2(a_4+a_5)u_{1,11} +
	(m_1+2m_2)\tau_{,1} |_{\Gamma} & = 0,
	\\
	\notag
	\sigma_{221}|_{\Gamma}  =
	a_1(u_{1,11}+u_{2,12}) +2a_3(u_{1,11} + u_{1,22}) +2a_4 u_{1,22} &
	\\ \label{NESBoundaryConditions2}
	+ 2a_5 u_{2,12} + m_1\tau_{,1} |_{\Gamma} & = 0,
	\\ \notag
	\sigma_{112}|_{\Gamma}  = \frac12
	a_1(u_{1,11}+u_{1,22}) + a_2(u_{1,11} + u_{2,21}) +2a_4 u_{2,12} & \\
	\label{NESBoundaryConditions3}   + a_5 ( u_{1,22} + u_{2,12}) +
	m_2\tau_{,1} |_{\Gamma} & = 0,
	\\
	\label{NESBoundaryConditions4}
	( \sigma_{112,1}, \sigma_{222,2} ) \cdot (\nu_1, \nu_2) |_{\Gamma}
	&= 0,
	\\
	\label{NESBoundaryConditions5}   (\tau_{,1}, \tau_{,2})^T \cdot
	(\nu_1, \nu_2)^T |_{\Gamma_2} = 0, \tau |_{\Gamma_1}  & = 0,
\end{align}
where $\nu = (\nu_1, \nu_2)^T$ is the unit exterior normal vector to $\Gamma = \partial \Omega$. 
To eliminate the trivial kernel, the ground space is then selected as the Hilbert space
\begin{equation}
	\spH_* =  \left( H^2(\Omega)/ \{1\} \right)^2 \times (L^2(\Omega)/
	\{1\})^2 \times \left(H^1(\Omega)/ \{1\}\right) \times
	\left(L^2(\Omega)/ \{1\}\right)
	\notag
\end{equation}
i.e., each component is taken as the orthogonal complement 
of the one-dimensional subspace of constant functions.
Without loss of generality, we equip this new space $\spH_*$ with the canonical unweighted inner product,
which is equivalent with our original weighted definition in Equation (\ref{EQUATION_SCALAR_PRODUCT_WEIGHTED}).

The operator $\opA$ is defined similar to (\ref{WPOpADefine}). 
The boundary conditions in Equations
(\ref{NESBoundaryConditions1})--(\ref{NESBoundaryConditions5}) are
incorporated into the weak definition of the domain $D(\opA)$ 
similar to Equation (\ref{WPOperatorDomain}).

\begin{theorem}
	Assume $a_3 + a_4 \ge 0$. 
	Then the system (\ref{NESElemSyst1})--(\ref{NESElemSyst3}) with boundary conditions
	(\ref{NESBoundaryConditions1})--(\ref{NESBoundaryConditions5}) is not exponentially stable.
\end{theorem}

\begin{proof}
We prove there exists a sequence $(\lambda_n)_n \subset \mathbb{R}$ with
\begin{equation}
	\lim_{n\to\infty} |\lambda_n| = \infty \notag
\end{equation}
as well as the sequences $(U_n)_n \subset D(\opA)$ and
$(F_n)_n\subset \spH_0$ such that
\begin{equation}
	(\ii \lambda_n - \opA)U_n = F_n 
	\text{ is  uniformly bounded w.r.t. }
	n\in\mathbb{N} \text{ and } \lim_{n\to\infty} \| U_n \|_{\spH_*} = \infty. \notag
\end{equation}
The lack of exponential stability will then follow from the well-known
Gearhart \& Pr\"uss' theorem (see e.g. \cite[Theorem 1.3.2, p. 4]{Liu1}).

For $\lambda\in\mathbb{R}$, the solution 
$U = (u_1, u_2, v_1, v_2, \tau, \theta)^T$ of the resolvent equation $\protect{(\lambda \ii - \opA)U = F}$ 
satisfies
\begin{align*}
	\lambda_n u_1 \ii - v_1   &= 0,
	\\
	\lambda_n u_2 \ii - v_2   &= 0,
	\\ 
	\lambda_n v_1 \ii - \frac{1}{\rho} \big((\lambda+2\mu)u_{1,11} +
	\mu u_{1,22} + (\lambda + \mu) u_{2,12} - \beta \dot{\tau}_{,1} -
	b_1 u_{1,1111} &
	\\ 
	- b_2 u_{1,1122} - b_3 u_{1,2222} - b_4 u_{2,1112} - b_5 u_{2,1222}&
	\\
	-(m_1+2m_2)(\tau_{,111} + \tau_{,122})\big)  &= f_3,
	\\
	\lambda_n v_2 \ii - \frac{1}{\rho} \big((\lambda + 2\mu)u_{2,22} +
	\mu u_{2,11} + (\lambda + \mu) u_{1,12} - \beta \dot{\tau}_{,2} - b_1 u_{2,2222} &
	\\
	- b_2 u_{2,1122} - b_3 u_{2,1111} - b_4
	u_{1,1112} - b_5 u_{1,1222} &
	\\
	-(m_1+2m_2)(\tau_{,112} + \tau_{,222})\big) &=  f_4,
	\\
	\lambda_n \tau \ii - \theta &= 0 ,
	\\
	\lambda_n \theta \ii - \frac{1}{a}\big((m_1 + 2m_2)(u_{1,111}+u_{1,122} + u_{2,112} +
	u_{2,222}) &  \\
	+ \tau_{,1} + \tau_{,2} - \beta (\dot{u}_{1,1} +
	\dot{u}_{2,2}) + \delta (\dot{\tau}_{,11} + \dot{\tau}_{,22}) \big) &=  f_6,
\end{align*}
where $f_3,f_4,f_6$ will be selected as
\begin{equation}
	f_3 = \sin (nx) \sin (ny), \quad
	f_4 = \cos (nx) \cos (ny), \quad
	f_6 = \cos (nx) \sin (ny)  \notag
\end{equation}
for $n \in \mathbb{N}$.
Eliminating $v_1, v_2, \theta$, we obtain for $u_1, u_2, \tau$ the following algebraic system
\begin{align}
	\notag
	-\lambda^2_n u_1 - \frac{1}{\rho} \big((\lambda+2\mu) u_{1,11} +
	\mu u_{1,22} + (\lambda + \mu) u_{2,12} &
	\\
	\notag
	-\beta \lambda_n \ii \tau_{,1} - b_1 u_{1,1111} - b_2 u_{1,1122} - b_3 u_{1,2222} 
	\\ 
	\label{NESReducedSystem1}
	- b_4 u_{2,1112} - b_5 u_{2,1222} -(m_1+2m_2)(\tau_{,111} + \tau_{,122})\big) &=  f_3,
	\\
	\notag
	-\lambda^2_n u_2  - \frac{1}{\rho} \big((\lambda + 2\mu) u_{2,22} +
	\mu u_{2,11} + (\lambda + \mu) u_{1,12} &
	\\
	\notag
	-\beta \lambda_n \ii \tau_{,2} - b_1 u_{2,2222}  - b_2 u_{2,1122} - b_3 u_{2,1111}
	\\ 
	\label{NESReducedSystem2}
	- b_4 u_{1,1112} - b_5 u_{1,1222} -(m_1+2m_2)(\tau_{,112} + \tau_{,222})\big) &=  f_4,
	\\
	\notag 
	-\lambda^2_n \theta \ii - \frac{1}{a} \big((m_1 + 2m_2)(u_{1,111}+u_{1,122} + u_{2,112} + u_{2,222}) &
	\\ 
	\label{NESReducedSystem3}
	+ \tau_{,1} + \tau_{,2} - \beta \lambda_n \ii (u_{1,1} +
	u_{2,2}) + \delta \lambda_n \ii (\tau_{,11} + \tau_{,22})\big) &=  f_6.
\end{align}
To solve Equations (\ref{NESReducedSystem1})--(\ref{NESReducedSystem3}), 
we employ the ansatz
\begin{equation}
	u_1 = A\sin (nx) \sin (ny), \quad
	u_2 = B\cos (nx) \cos (ny), \quad
	\tau = C\cos (nx) \sin (ny), \notag
\end{equation}
where $A, B, C$ will actually depend on $n$.
It should be pointed out that this choice is compatible with the
boundary conditions in Equations (\ref{NESBoundaryConditions1})--(\ref{NESBoundaryConditions5}).
Thus, system (\ref{NESReducedSystem1})--(\ref{NESReducedSystem3}) is
equivalent with finding $A, B, C$ such that
\begin{align}
	\notag 
	\big(-\lambda_n^2 \rho + (\lambda+ 3\mu) n^2 + 2b_2n^4\big) A - \big((\lambda +\mu) n^2 + 2b_4 n^4\big) B &
	\\ 
	\label{NESFinalSystem1}
	+ \big(-\beta \lambda_n n \ii + 2 (m_1+2m+2) n^3\big) C &=  \rho,
	\\
	\notag
	\big(-\lambda_n^2 \rho + (\lambda+ 3\mu) n^2 + 2b_2n^4\big) B - \big((\lambda +\mu) n^2 + 2b_4 n^4\big) A &
	\\ 
	\label{NESFinalSystem2}
	+ \big(\beta \lambda_n n \ii - 2(m_1 + 2m+2) n^3\big) C &=  \rho, \\
	\label{NESFinalSystem3}
	\big(2(m_1 + 2m_2) n^3 + \beta \lambda_n n \ii\big)(A - B) 
	+ (-\lambda_n^2 a + 2n^2 +2\lambda_n \delta n^2 \ii) C &=  a.
\end{align}
Let
\begin{equation}
	\lambda_n = \sqrt{\big(8(a_3+a_4) n^4 + 2\mu n^2 - 1\big)/\rho}. \notag
\end{equation}
It is easy to verify that the linear algebraic system
(\ref{NESFinalSystem1})--(\ref{NESFinalSystem3}) has a unique solution.
Summing up Equations (\ref{NESFinalSystem1}) and
(\ref{NESFinalSystem2}), we get
\begin{equation}
	(A + B) \big(-\lambda_n^2 \rho + 2\mu n^2 + 8(a_3 + a_4) n^4\big) = 2\rho, \notag
\end{equation}
where $a_3 + a_4 \ge 0$ follows from the positive definiteness of $C_{iIJKLj}$.
Equivalently,
\begin{equation}
	A + B = 2 \rho. \notag
\end{equation}
and, therefore, by virtue of Young's inequality,
\begin{align*} 
	\|U\|^2_{\spH_*} &\ge
	\|v_1\|^2_{L^2(\Omega)} + \|v_2\|^2_{L^2(\Omega)} 
	= \lambda_n^2 \left(\|u_1\|^2_{L^2(\Omega)} + \|u_2\|^2_{L^2(\Omega)}\right)
	\\
	&= \lambda_n^2 \frac{\pi^2}{2} (A^2 + B^2) 
	\ge \frac{\lambda_n^2 \pi^2}{2} \left(\frac{A+B}{2}\right)^2 
	= \frac{\lambda_n^2 \pi^2 \rho^2}{2} \to \infty,
\end{align*}
whereas we estimate
\begin{equation}
	\|F\|_{\spH^{\ast}}^{2} = \frac{3 \pi^{2}}{4} \text{ for any } n \in \mathbb{N}, \notag
\end{equation}
i.e., $(F_n)_{n}$ remains uniformly bounded w.r.t. $n \in \mathbb{N}$.
Hence, the claim follows.
\end{proof}

%%%%%%%%%%%%%%%%%%%%%%%%%%%%%%%%%%%%%%%%%%%%%%%%%%%%%%%%%%%%%%%%%%%%%%%%%%%%%%%%%%%%%%%%%%%%%%%%%%%%%%%%%%%%%%%%%%%%%%%%%%
%%%%%%%%%%%%%%%%%%%%%%%%%%%%%%%%%%%%%%%%%%%%%%%%%%%%%%%%%%%%%%%%%%%%%%%%%%%%%%%%%%%%%%%%%%%%%%%%%%%%%%%%%%%%%%%%%%%%%%%%%%

\section{Linear system with a frictional damping for the elastic part: Exponential stability}
\label{SECTION_EXPONENTIAL_STABILITY}
In Section \ref{SECTION_LACK_OF_EXPONENTIAL_STABILITY},
if no additional damping for the elastic component is present,
Equations (\ref{SystemTwoFrictionSecondCase1})--(\ref{SystemTwoFrictionSecondCase2})
were shown, in general, not to exhibit an exponential decay rate.
This justifies the necessity of introducing a damping mechanism for the elastic variable.
In the following, we consider a frictional damping
which leads to a system of partial differential equations reading as
\begin{align}
	\label{SystemTwoFrictionFirstCase1}
	\rho \ddot{u}_i &= \big(A_{iJKj} u_{j,K} - \beta_{Ji} \dot{\tau} -
	(C_{iJKLIj} u_{j,IL} + M_{iJKL} \tau_{,L} )_{,K}\big)_{,J} -
	E_{ij} \dot{u}_j,
	\\ 
	\label{SystemTwoFrictionFirstCase2}
	a \ddot{\tau} &=
	-\beta_{Ki} \dot{u}_{i,K} + m_{IJ}\dot{\tau}_{,IJ} + M_{jLKI} u_{j,LKI} +
	K_{IJ}\tau_{,IJ}
\end{align}
together with the boundary conditions (\ref{WPBoundaryConditions})
and initial conditions (\ref{WPInitialConditions}).
For this system, we assume
\begin{equation}
	E_{ij} \xi_{i} \xi_{j} \geq \alpha \xi_{i} \xi_{i} \text{ and }
	m_{IJ} \xi_{I} \xi_{J} \geq \alpha \xi_{I} \xi_{I} \text{ for some } \alpha > 0.
	\label{EQUATION_ASSUMPTION_POSITIVE_DEFINITENESS_OF_DAMPING}
\end{equation}

The natural first-order energy associated with the mild solution to
Equations (\ref{SystemTwoFrictionFirstCase1})--(\ref{SystemTwoFrictionFirstCase2}),
(\ref{WPBoundaryConditions})--(\ref{WPInitialConditions}) reads as
\begin{align}
	\label{ESFCEnergyDefinition}
	\begin{split}
		\mathcal{E}(t) &= \frac{1}{2} \int_{\Omega} (\rho \dot{u}_{i}^2 + a
		\dot{\tau}^2 + A_{iKLj} u_{i,K} u_{j,L} + C_{iIJKLj} u_{i,JI}
		u_{j,LK}) \, \mathrm{d}x
		\\ 
		&+ \frac12 \int_{\Omega}
		(K_{IJ} \tau_{,I} \tau_{,J} + 2 M_{iJKL} u_{i,KJ} \tau_{,L}) \, \mathrm{d}x,
	\end{split}
\end{align}
where, as before, the Einstein's summation convention is applied to terms like $\dot{u}_{i}^{2} = \dot{u}_{i} \dot{u}_{i}$, etc.
Due to the assumptions in Equation (\ref{EQUATION_ASSUMPTION_POSITIVE_DEFINITENESS_OF_DAMPING}),
there exists a number $c_{\mathcal{E}} > 0$ such that
\begin{equation}
	\label{ESFCEstimationForEnergy}
	\mathcal{E}(t) \ge c_{\mathcal{E}} \int_{\Omega} (\dot{u}_{i}^2 +
	\dot{\tau}^2 + u^2_{i,J} +  u^2_{i,JK} +
	\tau^2_{,I})(t,x) \, \mathrm{d}x.
\end{equation}

\begin{theorem}
	\label{ESFCTheorem}
	Let $U_0 \in \spH$. 
	There exist then positive constants $C$ and $c_0$ such that
	\begin{equation}
		\mathcal{E}(t) \le C \mathcal{E}(0) e^{-c_0 t}
		\text{ holds true for every } t \geq 0. \notag
	\end{equation}
\end{theorem}

\begin{proof}
Without loss of generality, let $U_0 \in D(\opA)$. 
Indeed, if this is not the case, 
due to the dense embedding of $D(\opA) \hookrightarrow \spH$, 
an appropriate approximating sequence from $D(\opA)$ can be selected. 
Further, let $(u, \dot{u}, \tau, \dot{\tau})^{T}$ denote the classical solution to
(\ref{SystemTwoFrictionFirstCase1})--(\ref{SystemTwoFrictionFirstCase2}),
(\ref{WPBoundaryConditions})--(\ref{WPInitialConditions}) for the initial data $U_0$. 
We construct a Lyapunov functional $\mathcal{F}$. 
Computing $\p \mathcal{E}(t)$, we get
\begin{align*}
	\p \mathcal{E}(t) &= \int_{\Omega} \left(-\frac{1}{\rho} E_{ij}
	\dot{u}_i \dot{u}_j + \frac{1}{a} m_{IJ} \dot{\tau}_{,IJ}
	\dot{\tau}\right) \, \mathrm{d}x
	\\
	&= -\int_{\Omega} \left(\frac{1}{\rho} E_{ij} \dot{u}_i \dot{u}_j
	+ \frac{1}{a} m_{IJ} \dot{\tau}_{,I} \dot{\tau}_{,J}\right) \
	\mathrm{d}x\le -\alpha \int_{\Omega}  \left(\dot{u}_i^2 +
	\dot{\tau}^2_{,I}\right) \, \mathrm{d}x.
\end{align*}
Letting
\begin{equation}
	\mathcal{F}_1(t) = \rho \int_{\Omega}  \dot{u}_i u_i \, \mathrm{d}x \text{ and }
	\mathcal{F}_2(t) = a \int_{\Omega} \dot{\tau} \tau \, \mathrm{d}x,
\end{equation}
after a partial integration, we arrive at
\begin{align*}
	\p \mathcal{F}_1(t) &= \int_{\Omega} \left(-A_{iKLj} u_{i,K} u_{j,L} -
	\beta_{Ji} \dot{\tau}_{,J} u_i - C_{iIJKLj} u_{i,JI}
	u_{j,LK}\right) \, \mathrm{d}x
	\\
	&+ \int_{\Omega }\left(-
	M_{iJKL} \tau_{, L} u_{i, KJ} - E_{ij} \dot{u}_i u_j + \rho
	\dot{u}^2\right) \, \mathrm{d}x,
	\\
	\p \mathcal{F}_2(t) = & \int_{\Omega}\left( \beta_{Ki} \dot{u}_{i}
	\tau_{,K} - m_{IJ} \dot{\tau}_{,I} \tau_{,J} - M_{iJKL} \tau_{, L}
	u_{i, KJ} - K_{IJ} \tau_{,I}\tau_{,J} + a \dot{\tau}^2\right) \, \mathrm{d}x.
\end{align*}
Let $R_0 > \max\{\beta_{Ji}, m_{IJ}, E_{ij}, a, \rho \}$. 
Utilizing the generalized Young's inequality and the first Poincar\'{e}'s inequality,  we obtain
\begin{align*}
	\p & (\mathcal{F}_1  + \mathcal{F}_2)(t) =  \int_{\Omega}\left(
	-A_{iKLj} u_{i,K} u_{j,L} - C_{iIJKLj} u_{i,JI} u_{j, LK} - 2M_{iJKL}
	\tau_{, L} u_{i, KJ}\right) \, \mathrm{d}x
	\\
	& + \int_{\Omega}  \left(\beta_{Ki} \dot{u}_{i} \tau_{,K} - K_{IJ} \tau_{,I}\tau_{,J} - \beta_{Ji} \dot{\tau}_{,J}
	u_{i} - m_{IJ} \dot{\tau}_{,I} \tau_{,J}  - E_{ij} \dot{u}_i u_j +
	\rho \dot{u}^2 + a \dot{\tau}^2\right) \, \mathrm{d}x
	\\
	\le & -\alpha \int_{\Omega} \left( u^2_{i,J} +  u^2_{i,KL} +
	\tau^2_{,I}\right) \ \mathrm{d}x + \frac{R_0}{2} \int_{\Omega}
	\left(\frac{1}{\varepsilon} \dot{u}^2_i + \varepsilon
	\tau^2_{,K}\right) + \left(\frac{1}{\varepsilon} \dot{\tau}^2_{,J} +
	\varepsilon u_i^2 \right) \, \mathrm{d}x
	\\
	& +\frac{R_0}{2}\int_{\Omega} \left(\frac{1}{\varepsilon} \dot{\tau}^2_{,I} + \varepsilon
	\tau_{,J}^2 \right) + \left(\frac{1}{\varepsilon} \dot{u}^2_{i} +
	\varepsilon u_{j}^2 \right) + \left(\dot{u}^2 + \dot{\tau}^2\right) \, \mathrm{d}x
	\\
	\le & -\alpha \int_{\Omega}\left( u^2_{i,J} + u^2_{i,KL} +
	\tau^2_{,I}\right) \ \mathrm{d}x + \frac{R_0n}{2} \int_{\Omega}
	\left(\frac{1}{\varepsilon} \dot{u}^2 + \varepsilon \tau^2_{,K} +
	\frac{1}{\varepsilon} \dot{\tau}^2_{,I}\right) \, \mathrm{d}x
	\\
	&+ \frac{R_0n}{2} \int_{\Omega}
	\left(\varepsilon c_{F} u_{i,J}^2  + \frac{1}{\varepsilon}
	\dot{\tau}^2_{,I} + \varepsilon \tau_{,J}^2  + \frac{1}{\varepsilon}
	\dot{u}^2 + \varepsilon c_{F} u_{j,I}^2 + \dot{u}^2 +
	c_{F}\dot{\tau}_{,I}^2\right) \, \mathrm{d}x.
\end{align*}
Now, selecting $\varepsilon>0$ such that
\begin{equation}
	\frac{\alpha}{2} > \frac{R_0 n}{2} \max \{ 2\varepsilon, 2\varepsilon c_{F}\}, \notag
\end{equation}
we find
\begin{align*}
	\p (\mathcal{F}_1 + \mathcal{F}_2)(t) &\le -\frac{\alpha}{2}
	\int_{\Omega}\left( u^2_{i,J} + u^2_{i,KL} +
	\tau^2_{,I}\right) \ \mathrm{d}x
	\\  
	&+ \frac{R_0n}{2} \int_{\Omega}
	\left(\frac{2}{\varepsilon}+1\right) \dot{u}^2 +
	\left(\frac{2}{\varepsilon}+c_{F}\right) \dot{\tau}^2_{,I} \, \mathrm{d}x.
\end{align*}
Next, we define
\begin{equation}
	\mathcal{F}(t) = \mathcal{F}_1(t) + \mathcal{F}_2(t) + N \mathcal{E}(t) \notag
\end{equation}
for some $N> \frac12 + \frac{R_0n}{2\alpha}\left(\frac{2}{\varepsilon} + \max\{1,c_{F}\} \right)$ to be fixed later. 
Then
\begin{equation}
	\p \mathcal{F}(t) \le -\frac{\alpha}{2} \int_{\Omega}
	\left(u^2_{i,J} + u^2_{i,KL} +  \tau^2_{,I} + \dot{u}^2 +
	\dot{\tau}^2_{,I}\right) \, \mathrm{d}x. \notag
\end{equation}
Using the first Poincar\'{e}'s inequality and Equation (\ref{ESFCEstimationForEnergy}), we obtain
\begin{equation}
	\p \mathcal{F}(t) \le -\widehat{C} \mathcal{E}(t). \notag
\end{equation}
Taking into account
\begin{equation}
	\big|(\mathcal{F}_1 + \mathcal{F}_2)(t)\big| \le \widetilde{C} \mathcal{E}(t), \notag
\end{equation}
we conclude that
\begin{equation}
	(N - \widetilde{C}) \mathcal{E}(t) \le \mathcal{F}(t) \le (N+\widetilde{C}) \mathcal{E}(t). \notag
\end{equation}
If necessary, $N$ is increased to make $N-\widetilde{C}$ positive.
Gronwall's inequality now yields
\begin{equation}
	\mathcal{E}(t) \le \frac{1}{N-\widetilde{C}} \mathcal{F}(t) \le
	\frac{1}{N-\widetilde{C}} \mathcal{E}(0)
	e^{-\widehat{C}/(N+\widetilde{C})} t = C \mathcal{E}(0) e^{-c_0t}. \notag
\end{equation}
This completes the proof.
\end{proof}

%%%%%%%%%%%%%%%%%%%%%%%%%%%%%%%%%%%%%%%%%%%%%%%%%%%%%%%%%%%%%%%%%%%%%%%%%%%%%%%%%%%%%%%%%%%%%%%%%%%%%%%%%%%%%%%%%%%%%%%%%%
%%%%%%%%%%%%%%%%%%%%%%%%%%%%%%%%%%%%%%%%%%%%%%%%%%%%%%%%%%%%%%%%%%%%%%%%%%%%%%%%%%%%%%%%%%%%%%%%%%%%%%%%%%%%%%%%%%%%%%%%%%

\section{Hyperbolized system with a frictional damping in the elastic part: Global existence and exponential stability}
\label{SECTION_HYPERBOLIZED_SYSTEM}
In this last section, we want to study the impact
of a nonlinear frictional damping on `the' equations of thermoelasticity for nonsimple materials.
In contrast to Section \ref{SECTION_EXPONENTIAL_STABILITY},
we replace Equation (\ref{SystemTwoFrictionSecondCase2}) with a Cattaneo-like hyperbolic relaxation
(see, e.g., \cite{TaZh1998}).
Thus, the resulting system becomes purely hyperbolic
and we can use the well-known technique due to Lasiecka and Tataru
\cite{Lasiecka1} to utilize an observability inequality for the linear system
to obtaine a uniform decay for the nonlinear one.

Letting $\mathbf{q} = (q_{i})$ such that $(m_{IJ} q_{I, J})$ represents the heat flux,
for a (small) relaxation parameter $\kappa > 0$,
we consider the following semilinear initial-boundary value problem
\begin{align}
	\label{NLSystemWithNonLinDamping1}
	\rho \ddot{u}_i &= \left( A_{iJKj} u_{j,K} - \beta_{Ji} \dot{\tau} -
	(C_{iJKLIj} u_{j,IL} + M_{iJKL} \tau_{,L} )_{,K} \right)_{,J} - E(|\dot{u}|) \dot{u}_i,
	\\  
	\label{NLSystemWithNonLinDamping2}
	a \ddot{\tau} &=
	-\beta_{Ki} \dot{u}_{i,K} + m_{IJ}q_{I,J} + M_{jLKI} u_{j,LKI} + K_{IJ}\tau_{,IJ},
	\\ 
	\label{NLSystemWithNonLinDamping3}
	\kappa \dot{q}_i &= \dot{\tau}_{,i} - q_{i}
\end{align}
in $\Omega \times (0, \infty)$ subject to the boundary conditions
\begin{align}
	\label{NLSystemWithNonLinDampingBoundCond}
	u_i = 0, \quad u_{i,J} = 0, \quad \tau = 0 \text{ in } \Gamma \times (0, \infty)
\end{align}
and the initial conditions
\begin{align}
	\label{NLSystemWithNonLinDampingInitCond}
	\begin{split}
		u_i(\cdot, 0) &= u_i^0, \quad \dot{u}_i(\cdot, 0) = \dot{u}_i^0, \\ 
		\tau (\cdot, 0) &= \tau^0, \quad \dot{\tau}(\cdot, 0) = \dot{\tau}^0, \quad q_i(\cdot,0) = q^0_i
		\text{ in } \Omega.
	\end{split}
\end{align}
Note that if $\kappa = 0$ and $E(\cdot)$ is linear,
Equations (\ref{NLSystemWithNonLinDamping1})--(\ref{NLSystemWithNonLinDamping3})
reduce to (\ref{SystemTwoFrictionSecondCase1})--(\ref{SystemTwoFrictionSecondCase2}).

\subsection{Preliminaries}
We assume the function $E \colon [0, \infty) \to \mathbb{R}$ satisfies the following conditions
\begin{itemize}
	\item $E(s) > 0$ for $s > 0$,
	
	\item $s \mapsto E(s)s$ is continuous and monotonically increasing,
	
	\item $\lim\limits_{s \searrow 0} E(s)s = 0$,
	
	\item There exist $m,M$ such that
	\begin{equation} \label{NLAssumptionForE4}
		0 < m \le E(s) \le M \text{ for } s > 1,
	\end{equation}
	
	\item The function $g_i$ given by $g_i(v) = E(|v|)v_i$ is uniformly Lipschitz-continuous,  i.e.,
	\begin{equation}\label{NLAssumptionForE5}
		|g_i(\phi) - g_i(\psi) | \le L |\phi - \psi|
	\end{equation}
	for all $\phi,\psi\in \mathbb{R}^d$.
\end{itemize}

For Equations (\ref{NLSystemWithNonLinDamping1})--(\ref{NLSystemWithNonLinDamping3}),
we assume the evolution is taking place on the Hilbert space
\begin{align*}
	\mathcal{H}_{\kappa} &= \Big\{U = (u, v, \tau, \theta, q)^{T} \in 
	\big(H^2_{0}(\Omega)\big)^d \times \big(L^2(\Omega)\big)^d \times
	H_0^1(\Omega) \times L^2(\Omega) \times (L^2(\Omega))^d\Big\}
	\\
	&= \spH \times (L^2(\Omega))^d
\end{align*}
equipped with the scalar product
\begin{align*}%\label{WPScalarProduct}
	\langle U , U^* \rangle_{\mathcal{H}_\kappa} = \big\langle
	(u,v,\tau,\theta)^T , (u^*,v^*,\tau^*,\theta^*)^T\big\rangle_{\mathcal{H}} 
	+ \int_{\Omega} \kappa m_{ij} q_i^* q_j \, \mathrm{d}x
\end{align*}
for $U = (u, v, \tau, \theta, q), U^* = (u^*, v^*, \tau^*, \theta^*, q^*)\in \spH_{\kappa}$.
(See Equation (\ref{EQUATION_SCALAR_PRODUCT_WEIGHTED}).)
As before, one can easily prove the inner product 
is equivalent with the usual product on the Hilbert space $\mathcal{H}_{\kappa}$.

We rewrite the problem
(\ref{NLSystemWithNonLinDamping1})--(\ref{NLSystemWithNonLinDampingInitCond}) in the abstract form:
\begin{equation}
	\label{NLOperatorCauchyProblem}
	\dot{U}(t) + \mathcal{K}\big(U(t)\big) = 0 \text{ for } t > 0, \quad U(0) = U_0.
\end{equation}
Here, $\mathcal{K} = - \mathcal{L} + \mathcal{N}$ is a nonlinear operator
with the domain $D(\mathcal{K}) := D(\mathcal{L})$, where
\begin{equation}\label{NLDefineL}
	\mathcal{L} U = 
	\left(
	\begin{array}{c}
		v_i \\
		\frac{1}{\rho} \left( A_{iJKj} u_{j,K} - \beta_{Ji} \theta -
		(C_{iJKLIj} u_{j,IL} + M_{iJKL} \tau_{,L} )_{,K} \right)_{,J} \\
		\theta \\
		\frac{1}{a}\left( -\beta_{Ki} v_{i,K} + m_{IJ}q_{I,J} + M_{jLKI} u_{j,LKI} + K_{IJ}\tau_{,IJ}\right) \\
		\frac{1}{\kappa} \theta_{,i} \\
	\end{array}
	\right)
\end{equation}
for $u \in D(\mathcal{L})$ and
\begin{equation}
	\label{NLDefineN}
	\mathcal{N}(U) =
	\Big(0, 
	\frac{1}{\rho} E(|v|) v_i,
	0,
	0,
	\frac{1}{\kappa} q_{i}\Big)^{T} \text{ for } u \in D(\mathcal{N})
\end{equation}
as well as
\begin{align}
	D(\mathcal{L}) &= \Big\{U \in \mathcal{H} \,|\, \phantom{-}v_i \in H^1_0(\Omega), \, \theta \in H^1_0(\Omega), \label{NLOperatorLDomain} \\
	&\phantom{= \phantom{-}\Big\{U \in \mathcal{H} \,|\,}
	\big(A_{iJKj} u_{j,K} - \beta_{Ji} \theta \notag \\
	&\phantom{= \Big\{U \in \mathcal{H} \,|\,} - (C_{iJKLIj} u_{j,IL} 
	+ M_{iJKL} \tau_{,L} )_{,K} \big)_{,J} \in L^2(\Omega) \notag \\
	&\phantom{= \phantom{-}\Big\{U \in \mathcal{H} \,|\,}
	(m_{JI}q_{J} + M_{jLKI} u_{j,LK} + K_{IJ}\tau_{,J})_{,I} \in L^2(\Omega)\Big\}, \notag \\
	D(\mathcal{N}) &= \mathcal{H}_{\kappa}. \label{NLOperatorNDomain}
\end{align}

\begin{definition}
	Let $U_{0} \in \mathcal{H}$
	and assume $\mathcal{L}$ generates a $C_{0}$-semigroup $\big(S(t)\big)_{t \geq 0}$ on $\mathcal{H}_{\kappa}$.
	A function
	$U \in C^{0}\big([0, \infty), \mathcal{H}_{\kappa}\big)$
	satisfying the integral equation
	\begin{equation}
		U(t) = S(t) U_{0} + \int_{0}^{t} S(t - s) \mathcal{N}\big(U(s)\big) \mathrm{d}s \notag
	\end{equation}
	is referred to as a mild solution to Equation (\ref{NLOperatorCauchyProblem}).
	If $U$ additionally satisfies
	\begin{equation}
		U \in H^{1}_{\mathrm{loc}}(0, \infty; \mathcal{H}_{\kappa}) \cap
		L^{2}_{\mathrm{loc}}\big(0, \infty; D(\mathcal{L})\big), \notag
	\end{equation}
	we call $U$ a strong solution.
\end{definition}

By $\nlE(t)$ we denote the corresponding energy
\begin{align*}%\label{NLEnergyDefinition}
	\nlE(t) &= \frac12 \int_{\Omega} \left(\rho |\dot{u}|^2 + a
	\dot{\tau}^2 + A_{iKLj} u_{i,K} u_{j,L} + C_{iIJKLj} u_{i,JI}
	u_{j,LK} + K_{IJ}
	\tau_{,I} \tau_{,J}\right) \, \mathrm{d}x \\
	&+ \frac12 \int_{\Omega}\left( 2 M_{iJKL} u_{i,KJ} \tau_{,L} +
	\kappa |q|^2\right) \, \mathrm{d}x
\end{align*}
associated with a mild solution $U$ to nonlinear Equation (\ref{NLOperatorCauchyProblem}).
We observe
\begin{align}
	\langle \mathcal{L}V , V \rangle_{\mathcal{H}_\kappa} &= 0 
	\text{ for any } V \in D(\mathcal{L}) \text{ and } 
	\label{NLDissipOfNonLinSyst1} \\
	\p \nlE(t) &= -\int_{\Omega} E(|v|) |v|^2  + m_{ij} q_i q_j \, \mathrm{d}x
	\text{ a.e. in } (0, \infty)
	\label{NLDissipOfNonLinSyst2}
\end{align}
for any strong solution $U$ to Equation (\ref{NLOperatorCauchyProblem}).

\begin{remark}
	By virtue of Stone's theorem
	(cf. \cite{St1932}, \cite[Theorem 3.8.6, p. 105]{TuWei2009}),
	Equation (\ref{NLDissipOfNonLinSyst1}) is equivalent with $\mathcal{A}$ being skew-adjoint.
\end{remark}

A straighforward adaption of the proof of Lemma \ref{WPMainLemma} yields:
\begin{lemma}
	\label{NLMainLemma}
	The operator $\mathcal{L}$ defined in Equation (\ref{NLDefineL}) is the infinitesimal
	generator of a $C_0-$semigroup of contractions on $\mathcal{H}_\kappa$.
\end{lemma}

For the nonlinear Cauchy problem (\ref{NLOperatorCauchyProblem}),
we have the following existence and uniqueness theorem.
\begin{theorem}
	For $U_0 \in \spH_\kappa$, 
	there exist unique mild solution $U$ to Equation (\ref{NLOperatorCauchyProblem}).
	If $U_{0} \in D(\opK)$, the mild solution is strong.
\end{theorem}

\begin{proof}
	By virtue of Equation (\ref{NLAssumptionForE5}),
	\begin{align*}
		\big\|\mathcal{N}(U^{1}) - \mathcal{N}(U^{2})\big\|_{\mathcal{H}_{\kappa}} &=
		\Big(\rho \int_{\Omega} \big|E(|v^{1}|) v^{1} - E(|v^{2}|) v^{2}\big|^{2} \, \mathrm{d}x\Big)^{1/2} \\
		&\leq L \sqrt{\rho} \|v^{1} - v^{2}\|_{L^{2}}
		\leq L \big\|U^{1} - U^{2}\|_{\mathcal{H}_{\kappa}} \notag
	\end{align*}
	for $U^{i} = (u^{i}, v^{i}, \tau^{i}, \theta^{i}, q^{i})^{T}$, $i = 1, 2$,
	i.e., the mapping $\mathcal{N}$ is globally Lipschitzian on $\mathcal{H}_{\kappa}$.
	Hence, the first claim is a direct consequence of \cite[Theorem 1.2, p. 184]{Pazy1}.
	Since $\mathcal{H}_{\kappa}$, being a Hilbert space, is reflexive,
	the second claim readily follows from \cite[Theorem 1.4, p. 189]{Pazy1}.
\end{proof}

By utilizing Lemma \ref{NLMainLemma}
and constructing a Lyapunov's functional similar to that one in Section \ref{SECTION_EXPONENTIAL_STABILITY} (see also \cite{Po2014}),
we further obtain the following linear stability theorem.
\begin{theorem}
	\label{THEOREM_LINEAR_HYPERBOLIZED_SYSTEM_EXPONENTIAL_STABILITY}
	For
	\begin{equation*}
		\label{NLDefineD}
		\mathcal{D} \colon \mathcal{H}_{\kappa} \to \mathcal{H}_{\kappa}, \quad
		\mathcal{D} U := \Big(
		0,
		\frac{1}{\rho} E_{ij} v_j,
		0,
		0,
		\frac{1}{\kappa} q_{i}
		\Big)^{T} \text{ for } U \in \mathcal{H}_{\kappa},
	\end{equation*}
	the unique (mild or classical) solution to the Cauchy problem
	\begin{equation}
		\dot{U}(t) = \opL U(t) - \mathcal{D}U(t) \text{ for } t > 0, \quad
		U(0) = U_0
	\end{equation}
	is exponentially stable.
\end{theorem}

\subsection{Uniform stability}
Recalling Equation (\ref{NLDissipOfNonLinSyst1})
and using Theorem \ref{THEOREM_LINEAR_HYPERBOLIZED_SYSTEM_EXPONENTIAL_STABILITY},
we apply Theorem \ref{APPTheorem} in Appendix \ref{SECTION_APPENDIX} to obtain 
the following observability result for the linear part of Equation (\ref{NLOperatorCauchyProblem}).
\begin{theorem}
	\label{ObservInequality}
	There exist a time period $T_{0} > 0$ and a positive constant $C$ such that
	for any $T \geq T_{0}$ 
	every mild solution $(u, \dot{u}, \tau, \dot{\tau}, q)^{T}$ to the Cauchy problem
	\begin{equation}
		\dot{U}(t) = \opL U(t) \text{ for } t > 0, \quad
		U(0) = U_{0} \notag
	\end{equation}
	satisfies
	\begin{equation}
		\mathcal{E} (0) \le C \int_0^T 
		\int_{\Omega} \big(|\dot{u}|^2 + |q|^2\big) \, \mathrm{d}x \mathrm{d}t. \notag
	\end{equation}
\end{theorem}

In the following, we adapt the techniques introduced by Lasiecka and Tataru \cite{Lasiecka1}.
The technical difficulties in our case are due to the big system size
and mixed-order structure of Equations (\ref{NLSystemWithNonLinDamping1})--(\ref{NLSystemWithNonLinDamping3}).
\begin{lemma}
	\label{NLLemaEstimationForEnergy}
	For any $T > 0$ and $K > 0$, 
	there exists a positive constant $C$ such that
	\begin{equation}
		\nlE(T) \le C \int_0^T \int_{\Omega}\left( E^2(|\dot{u}|)|\dot{u}|^2 
		+ |\dot{u}|^2 + |q|^2 \right) \, \mathrm{d}x \mathrm{d}t \notag
	\end{equation}
	for any strong solution $\Psi = (u_i, \dot{u}_i, \tau, \dot{\tau}, q)^{T}$
	to Equation (\ref{NLOperatorCauchyProblem}) additionally satisfying $\nlE(0) \le K$.
\end{lemma}

\begin{proof}
	Suppose the converse is true which means there exist $T > 0$ and $K > 0$ together with
	a sequence of strong solutions 
	$\big((u^{(m)}, \dot{u}^{(m)}, \tau^{(m)}, \dot{\tau}^{(m)}, q^{(m)})\big)_{m}^{T}$ such that $\nlE_m(0)\le K$
	to the Cauchy problem (\ref{NLOperatorCauchyProblem}), 
	whereas
	\begin{equation}
		\frac{\nlE_m(T)}{\int_0^T \int_{\Omega}
		E^2(|\dot{u}^{(m)}|)|\dot{u}^{(m)}|^2 + |\dot{u}^{(m)}|^2 +
		|q^{(m)}|^2 \, \mathrm{d}x \mathrm{d}t} \to \infty \text{ as } m \to \infty. \notag
	\end{equation}
	Denote
	\begin{equation}
		\Pi_m(\cdot) := 
		\int_0^T \int_{\Omega} E^2(|\dot{u}^{(m)}|)|\dot{u}^{(m)}|^2 + |\dot{u}^{(m)}|^2 +
		|q^{(m)}|^2 \, \mathrm{d}x \mathrm{d}t. \notag
	\end{equation}
	Using Equation (\ref{NLDissipOfNonLinSyst2}), we then get 
	\begin{equation*}
		\frac{K}{\Pi_m} \ge \frac{\nlE_m(0)}{\Pi_m} \ge
		\frac{\nlE_m(T)}{\Pi_m} \to \infty \text{ as } m \to \infty,
	\end{equation*}
	and, thus, $\Pi_m \to 0$, i.e.,
	\begin{equation*}
		\int_{0}^{T} \|\dot{u}^{(m)}\|_{L^2}^{2} \mathrm{d}t \to 0 \text{ and } 
		\int_{0}^{T} \|q^{(m)}\|_{L^2}^{2} \mathrm{d}t \to 0 \text{ as } m \to \infty.
	\end{equation*}
	Letting 
	\begin{equation}
		\nu_m = \sqrt{\nlE_m(0)}, \quad 
		\bar{u}^{(m)} = \frac{u^{(m)}}{\nu_m}, \quad
		\bar{\tau}^{(m)} = \frac{\tau^{(m)}}{\nu_m}, \quad
		\bar{q}^{(m)} = \frac{q^{(m)}}{\nu_m}, \notag
	\end{equation}
	we observe $\big(\bar{u}^{(m)}, \bar{\tau}^{(m)}, \bar{q}^{(m)}\big)^{T}$ is a strong solution to
	\begin{align}
		\label{NLSystemNormalized1}
		\rho \ddot{u}_i &= \big( A_{iJKj} u_{j,K} - \beta_{Ji} \dot{\tau} -
		(C_{iJKLIj} u_{j,IL} + M_{iJKL} \tau_{,L} )_{,K}\big)_{,J} -
		E(|\dot{u}|\nu_m) \dot{u}_i
		\\ 
		\label{NLSystemNormalized2}
		a \ddot{\tau} &=
		-\beta_{Ki} \dot{u}_{i,K} + m_{IJ}q_{I,J} + M_{jLKI} u_{j,LKI} +
		K_{IJ}\tau_{,IJ},
		\\ 
		\label{NLSystemNormalized3}
		\kappa \dot{q}_i &= \dot{\tau}_{,i} - q_{i}
	\end{align}
	in $\Omega \times (0, \infty)$ subject to the boundary conditions
	\begin{align*}
		u_i = 0, \quad u_{i,J} = 0, \quad \tau = 0 \text{ on } \Gamma \times (0, \infty)
	\end{align*}
	and the initial conditions
	\begin{align}
		u_i = \frac{u_i^0}{\nu_{m}}, \quad 
		\dot{u}_i = \frac{\dot{u}_i^0}{\nu_{m}}, \quad
		\tau = \frac{\tau^0}{\nu_{m}}, \quad 
		\dot{\tau} = \frac{\dot{\tau}^0}{\nu_{m}}, \quad
		q_i = \frac{q^0_i}{\nu_{m}} \text{ in } \Omega \times \{0\}. \notag
	\end{align}
	By $\normE_m(t)$ denote the energy of the solution to system 
	(\ref{NLSystemNormalized1})--(\ref{NLSystemNormalized3}). 
	Note that
	\begin{equation}
		\label{NLInitEnergyOfNonLinSys}
		\normE_m(0) = \frac{\nlE_m(0)}{\nu_m} = 1.
	\end{equation}
	Now,
	\begin{align*}
		\int_0^T \int_{\Omega} \big(E^2(|\dot{u}^{(m)}|)|\dot{\bar{u}}^{(m)}|^2 +
		|\dot{\bar{u}}^{(m)}|^2 + |\bar{q}^{(m)}|^2 \big) \mathrm{d}x \mathrm{d}t =
		\frac{\Pi_m}{\nu^2_m} = \frac{\Pi_m}{\nlE_m(0)} \to 0
	\end{align*}
	as $m \to \infty$. Therefore,
	\begin{equation}
		\label{NLTauToZero}
		\int_{0}^{T} \|\dot{\bar{u}}^{(m)}\|_{L^2}^{2} \mathrm{d}t \to 0 \text{ and } 
		\int_{0}^{T} \|\bar{q}^{(m)}\|_{L^2}^{2} \mathrm{d}t
		\to 0 \text{ as } m \to \infty
	\end{equation}
	and also
	\begin{equation}\label{NLDampingToZero}
		\int_{0}^{T} \big\|E\big(|\dot{u}^{(m)}|\big) |\dot{\bar{u}}^{(m)}|\big\|_{L^2}^{2} \mathrm{d}t \to 0 
		\text{ as } m \to \infty.
	\end{equation}
	
	Now, consider the following linear system
	\begin{align}
		\label{NLCorrespLinSystem1}
		\rho \ddot{v}_i &= \big(A_{iJKj} v_{j,K} - \beta_{Ji} \dot{\omega} -
		(C_{iJKLIj} v_{j,IL} + M_{iJKL} \omega_{,L} )_{,K}\big)_{,J},
		\\ 
		\label{NLCorrespLinSystem2}
		a \ddot{\omega} &=
		-\beta_{Ki} \dot{v}_{i,K} + m_{IJ}r_{I,J} + M_{jLKI} v_{j,LKI} +
		K_{IJ}\omega_{,IJ},
		\\
		\label{NLCorrespLinSystem3} 
		\kappa \dot{r}_i &= \dot{\omega}_{,i} - r_{i}
	\end{align}
	subject to the boundary conditions
	\begin{equation}
		\label{NLCorrespLinSystemBoundCond}
		v_i = 0, \quad \ v_{i,J} = 0, \quad \omega = 0 \text{ on } \Gamma \times (0, \infty)
	\end{equation}
	and the initial conditions
	\begin{align}
		v_i = \frac{u_i^0}{\nu_{m}}, \quad 
		\dot{v}_i = \frac{\dot{u}_i^0}{\nu_{m}}, \quad
		\omega = \frac{\tau^0}{\nu_{m}}, \quad 
		\dot{\omega} = \frac{\dot{\tau}^0}{\nu_{m}}, \quad
		r_i = \frac{q^0_i}{\nu_{m}} \text{ in } \Omega \times \{0\}.
		\label{NLCorrespLinSystemInitCond}
	\end{align}
	Letting
	$\big(\bar{v}^{(m)}, \bar{\omega}^{(m)}, \bar{r}^{(m)}\big)^{T}$ 
	denote the strong solution to Equations
	(\ref{NLCorrespLinSystem1})--(\ref{NLCorrespLinSystemInitCond}) and
	$\lE_m(t)$ be the corresponding energy,
	the initial conditions (\ref{NLCorrespLinSystemInitCond}) imply
	\begin{equation}
		\label{NLEqualityOfEnergies}
		\lE_m(0) = \normE_m(0). \notag
	\end{equation}
	Exploring Equations (\ref{NLSystemNormalized1})--(\ref{NLSystemNormalized3})
	and (\ref{NLCorrespLinSystem1})--(\ref{NLCorrespLinSystemInitCond}),
	we deduce that 
	$w^{(m)}_i = \bar{u}^{(m)}_i - \bar{v}^{(m)}_i$,
	$\chi^{(m)} = \bar{\tau}^{(m)} - \bar{\omega}^{(m)}$,
	$p^{(m)}_i = \bar{q}^{(m)}_i - \bar{r}^{(m)}_i$ 
	strongly solves the `incremental' system
	\begin{align}
		\label{NLMixedSystem1}
		\rho \ddot{w}_i &= \big(A_{iJKj} w_{j,K} - \beta_{Ji} \dot{\chi} -
		(C_{iJKLIj} w_{j,IL} + M_{iJKL} \chi_{,L} )_{,K}\big)_{,J}
		\\ 
		\notag &- E(|\dot{\bar{u}}^{(m)}\nu_m|) \dot{\bar{u}}_i^{(m)}
		\\
		\label{NLMixedSystem2}
		a \ddot{\chi} &=
		-\beta_{Ki} \dot{w}_{i,K} + m_{IJ}p_{I,J} + M_{jLKI} w_{j,LKI} +
		K_{IJ}\chi_{,IJ},
		\\
		\label{NLMixedSystem3}
		\kappa \dot{p}_i &= \dot{\chi}_{,i} - p_{i},
	\end{align}
	together with the boundary conditions
	\begin{equation}
		\label{NLMixedSystemBoundVal}
		w_i = 0, \quad w_{i,J} = 0, \quad \chi = 0 \text{ on } \Gamma \times (0, \infty)
	\end{equation}
	and the initial conditions
	\begin{equation}
		\label{NLMixedSystemInitVal}
		w_i(\cdot, 0) = \dot{w}_i(\cdot, 0) = 0, \quad
		\chi(\cdot, 0) = \dot{\chi}(\cdot,0) = 0, \quad
		p_i(\cdot, 0) = 0 \text{ in } \Omega.
	\end{equation}
	With $\mE_m(t)$ denoting the associated energy,
	Equation (\ref{NLMixedSystemInitVal}) implies $\mE_m(0) = 0$,
	for which we easily verify
	\begin{equation}
		\p \mE_m(t) 
		= \int_{\Omega} \Big(-E(|\dot{u}^{(m)}|)
		\dot{\bar{u}}_i^{(m)} \dot{w}_i^{(m)} - m_{ij} \bar{q}^{(m)}_j p_i\Big) \, \mathrm{d}x.
		\notag
	\end{equation}
	Therefore,
	\begin{equation}
		\begin{split}
			\mE_m(t) &= \mE_m(t) - \mE_m(0) \\
			&= -\int_0^t \int_{\Omega} \left(
			E(|\dot{u}^{(m)}|) \dot{\bar{u}}_i^{(m)} \dot{w}_i^{(m)} + m_{IJ}
			\bar{q}^{(m)}_j p_i \right) \, \mathrm{d}x \mathrm{d}t.
		\end{split}
		\label{EQUATION_DEFINITION_ENERRGY_DIFF}
	\end{equation}
	Note that $\int_{0}^{T} \|\dot{w}^{(m)}\|_{L^2}^{2} \mathrm{d}t$ and 
	$\int_{0}^{T} \|p^{(m)}\|_{L^2}^{2} \mathrm{d}t$ are bounded w.r.t. $m \in \mathbb{N}$. 
	Indeed, a straighforward computation yields
	\begin{equation}
		\p \lE_m(t) = 0. \notag
	\end{equation}
	Thus, $\lE_m(t) = \lE_m(0)$ and, therefore,
	\begin{equation}
		\frac{\rho}{2} \int_{\Omega} |\dot{\bar{v}}^{(m)}|^2 \, \mathrm{d}x +
		\frac{\kappa\alpha}{2} \int_{\Omega} |\bar{r}^{(m)}|^2 \, \mathrm{d}x
		\le \lE_m(t) = \lE_m(0) = \normE_m(0) = 1, \notag
	\end{equation}
	whence
	\begin{equation}
		\int_0^T \int_{\Omega} \Big(\frac{\rho}{2}
		|\dot{\bar{v}}^{(m)}|^2 + \frac{\kappa\alpha}{2} |\bar{r}^{(m)}|^2\Big) \mathrm{d}x \mathrm{d}t \le T.
		\notag
	\end{equation}
	Estimating
	\begin{align*}
		\Big(\int_{0}^{T} \|\dot{w}^{(m)}_i\|_{L^2}^{2} \mathrm{d}t\Big)^{1/2} &\le 
		\Big(\int_{0}^{T} \|\dot{\bar{u}}^{(m)}_i\|_{L^2}^{2} \mathrm{d}t\Big)^{1/2} + 
		\Big(\int_{0}^{T} \|\dot{\bar{v}}^{(m)}_i\|_{L^2}^{2} \mathrm{d}t\Big)^{1/2}
		\text{ and } \\
		\Big(\int_{0}^{T} \|p^{(m)}_i\|_{L^2}^{2} \mathrm{d}t\Big)^{1/2} &\le 
		\Big(\int_{0}^{T} \|\bar{q}^{(m)}_i\|_{L^2}^{2} \mathrm{d}t\Big)^{1/2} + 
		\Big(\int_{0}^{T} \|\bar{r}^{(m)}_i\|_{L^2}^{2} \mathrm{d}t\Big)^{1/2}
	\end{align*}
	and recalling Equation (\ref{NLTauToZero}), we arrive at the boundness of 
	$\int_{0}^{T} \|\dot{w}^{(m)}\|_{L^2}^{2} \mathrm{d}t$ and 
	$\int_{0}^{T} \|p^{(m)} \|_{L^2}^{2} \mathrm{d}t$ w.r.t. $m \in \mathbb{N}$.
	Using this fact and applying Cauchy \& Schwarz' inequality to Equation (\ref{EQUATION_DEFINITION_ENERRGY_DIFF}),
	we further find a positive constant $C$ such that
	\begin{align*}
		\mE_m(t) &\le \Big(\int_0^t\int_{\Omega} E^2(|\dot{u}^{(m)}|)
		|\dot{\bar{u}}_i^{(m)}|^2 \, \mathrm{d}x\mathrm{d}t\Big)^{1/2}
		\Big(\int_0^t\int_{\Omega}
		|\dot{w}_i^{(m)}|^2 \, \mathrm{d}x \mathrm{d}t\Big)^{1/2}
		\\
		&\phantom{\leq} + \max_{i,j} |m_{ij}| \Big(\int_0^t\int_{\Omega} |\bar{q}^{(m)}_j|^2
		\, \mathrm{d}x\mathrm{d}t\Big)^{1/2}
		\Big(\int_0^t\int_{\Omega} |p_i|^2 \ \mathrm{d}x\mathrm{d}t
		\Big)^{1/2}
		\\
		&\le C \Big(\int_0^T\int_{\Omega} E^2(|\dot{u}^{(m)}|)
		|\dot{\bar{u}}_i^{(m)}|^2 \, \mathrm{d}x \mathrm{d}t\Big)^{1/2}  
		+ C
		\Big(\int_0^T\int_{\Omega} |\bar{q}^{(m)}_j|^2 \, \mathrm{d}x \mathrm{d}t\Big)^{1/2}
		\\
		&\le C \left(\Big(\int_{0}^{T} \|E(|\dot{u}^{(m)}|) |\dot{\bar{u}}^{(m)}|\|_{L^2}^{2} \mathrm{d}t\Big)^{1/2}
		+ \Big(\int_{0}^{T}\|\bar{q}^{(m)}\|_{L^2}^{2} \mathrm{d}t\Big)^{1/2}\right).
	\end{align*}
	Hence,
	\begin{equation}
		\label{NLMixedEnergyToZero}
		\max_{t \in [0, T]} \mE_m(t) \le
		C \bigg(\Big(\int_{0}^{T} \|E(|\dot{u}^{(m)}|) |\dot{\bar{u}}^{(m)}|\|_{L^2}^{2} \mathrm{d}t\Big)^{1/2}
		+\Big(\int_{0}^{T}\|\bar{q}^{(m)}\|_{L^2}^{2} \mathrm{d}t\Big)^{1/2}\bigg).
	\end{equation}
	Further, using the trivial inequality 
	$a^2 = (a-b + b)^2 \le 2((a-b)^2 + b^2)$, we get
	\begin{align*}
		\int_{\Omega} \big(|\dot{\bar{v}}^{(m)}|^2 + |\bar{r}^{(m)}|^2 \big) \mathrm{d}x 
		&\le 2\int_{\Omega} \big(|\dot{\bar{u}}^{(m)} -
		\dot{\bar{v}}^{(m)}|^2 + |\bar{q}^{(m)} - \bar{r}^{(m)}|^2 +
		|\dot{\bar{u}}^{(m)}|^2 + |\bar{q}^{(m)}|^2\big) \mathrm{d}x
		\\
		&= 2 \int_{\Omega} \big(|\dot{w}^{(m)}|^2 + |p^{(m)}|^2 +
		|\dot{\bar{u}}^{(m)}|^2 + |\bar{q}^{(m)}|^2\big) \mathrm{d}x
		\\
		&\le C_1 \mE_m(t) + 2
		\int_{\Omega} \big(|\dot{\bar{u}}^{(m)}|^2 + |\bar{q}^{(m)}|^2\big) \mathrm{d}x
		\\
		&\le C_1 \max_{t \in [0, T]} \mE_m(t) + 
		2 \int_{\Omega} \big(|\dot{\bar{u}}^{(m)}|^2 + |\bar{q}^{(m)}|^2\big) \mathrm{d}x.
	\end{align*}
	Integrating the latter inequality w.r.t. $t$ and 
	using Equations (\ref{NLTauToZero}), (\ref{NLDampingToZero}) and (\ref{NLMixedEnergyToZero}), we obtain
	\begin{equation}
		\int_0^T \int_{\Omega} \big(|\dot{\bar{v}}^{(m)}|^2 + |\bar{r}^{(m)}|^2\big)
		\mathrm{d}x \mathrm{d}t \to 0 \text{ as } m \to \infty. \notag
	\end{equation}
	On the other hand,  recalling Equations
	(\ref{NLInitEnergyOfNonLinSys}), (\ref{NLEqualityOfEnergies}) and Theorem \ref{ObservInequality},
	we find
	\begin{equation}
		1 = \normE_m(0) = \lE_m(0) \le c \int_{\Omega}
		\big(|\dot{\bar{v}}^{(m)}|^2 + |\bar{r}^{(m)}|^2\big) \mathrm{d}x \mathrm{d}t \to 0
		\text{ as } m \to \infty, \notag
	\end{equation}
	which is a contradiction to our original assumption.
\end{proof}

To proceed further, let us introduce some notations.
Thanks to our assumptions on the function $E$,
according to \cite{Lasiecka1}, there exists a real-valued function $h \colon [0, \infty) \to [0, \infty)$,
being concave, strictly increasing and satisfying $h(0) = 0$ and
\begin{equation}
	\label{NLDefinOfH}
	h(s^2E(s)) \ge s^2 + E^2(s)s^2 \text{ for } s \in [0, 1].
\end{equation}
Further, we define a function $r$ by means of
\begin{equation}
	\label{NLDefinOfR}
	r(s) = h\left(\frac{s}{|Q|}\right) \text{ for } s \ge 0,
\end{equation}
where $Q = \Omega \times (0, T)$ and $|Q|$ is the standard Borel measure of $Q$.
Further, let
\begin{equation}
	\label{NLDefinOfP}
	p(s) = (cI+r)^{-1} \left(Ms\right) \text{ for } s \ge 0,
\end{equation}
where $M$ and $c$ are some positive constants to be defined later.
Finally, let
\begin{equation}
	\label{NLDefinOfQ}
	q(s) = s - (I+p)^{-1} \left(s\right) \text{ for } s \ge 0.
\end{equation}

We quote the following lemma due to Lasiecka and Tataru \cite{Lasiecka1} we use below.
\begin{lemma}
	\label{LemmaLasiecka}
	Let the functions $p, q$ be defined as above. 
	For any number sequence $(s_n)_n \subset (0, \infty)$ with $s_{m+1} + p(s_{m+1}) \le s_m$,
	we have
	\begin{equation}
		s_m \le S(m) \text{ for every } m\in \mathbb{N}, \notag
	\end{equation}
	where $S(t)$ is a solution of the scalar Cauchy problem
	\begin{equation}
		\label{NLDefineS}
		\dot{S}(t) + q\big(S(t)\big) = 0 \text{ for } t > 0, \quad
		S(0) = s_0.
	\end{equation}
	Moreover, if $p$ satisfies $p(s) > 0$ for $s > 0$, then
	\begin{equation}
		\lim_{t \to \infty} S(t) = 0. \notag
	\end{equation}
\end{lemma}

Now, we are in position to prove the uniform stability result for the Cauchy problem (\ref{NLOperatorCauchyProblem}).
\begin{theorem}
	\label{NLTheoremOnStability}
	For a sufficiently large number $T > 0$, 
	we have
	\begin{equation}
		\nlE(t) \le S\left(\frac{t}{T}-1\right) \text{ for } t > T. \notag
	\end{equation}
	Moreover,
	$S(t) \to 0$ as $t \to \infty$, where $S(t)$ solves the ODE (\ref{NLDefineS}).
\end{theorem}

\begin{proof}
	Consider the sets
	\begin{align*}
		Q_1 &= \big\{(x, t) \in \Omega \times (0,T) \,|\, \big|\dot{u}(x, t)\big| > 1 \big\},
		\\
		Q_2 &= \big\{(x, t) \in \Omega \times (0,T) \,|\, \big|\dot{u}(x, t)\big| \le 1 \big\}.
	\end{align*}
	Due to the additivity property of Lebesgue integral, we have
	\begin{align*}
		\int_0^T \int_{\Omega} \big(E^2(|\dot{u}|)|\dot{u}|^2 + |\dot{u}|^2 \big) \mathrm{d}x \mathrm{d}t 
		&= \int_{Q_1} \big(E^2(|\dot{u}|)|\dot{u}|^2 + |\dot{u}|^2\big) \mathrm{d}x \mathrm{d}t
		\\
		&+\int_{Q_2} \big(E^2(|\dot{u}|)|\dot{u}|^2 + |\dot{u}|^2\big) \mathrm{d}x \mathrm{d}t.
	\end{align*}
	On the set with $|\dot{u}| > 1$, Equation (\ref{NLAssumptionForE4}) implies
	\begin{align*}
		E^2(|\dot{u}|)|\dot{u}|^2 \le M E(|\dot{u}|)|\dot{u}|^2 \text{ and }
		|\dot{u}|^2 \le \frac{1}{m} E(|\dot{u}|)|\dot{u}|^2
	\end{align*}
	with $m, M$ from Equation (\ref{NLAssumptionForE4}).
	Therefore,
	\begin{align}
		\notag
		\int_{Q_1} \big(E^2(|\dot{u}|)|\dot{u}|^2 + |\dot{u}|^2\big) \mathrm{d}x \mathrm{d}t
		&\le \left(M + \frac{1}{m}\right) \int_{Q_1} E(|\dot{u}|)|\dot{u}|^2 \, \mathrm{d}x \mathrm{d}t
		\\ 
		\label{NLIneqForQ1}
		&\le c_1 \int_{Q} E(|\dot{u}|)|\dot{u}|^2 \, \mathrm{d}x \mathrm{d}t,
	\end{align}
	where $c_1 = M + m^{-1}$.

	Next, consider the integral over $Q_2$. 
	Using Equation (\ref{NLDefinOfH}), we estimate
	\begin{align*}
		\int_{Q_2} E^2(|\dot{u}|)|\dot{u}|^2 + |\dot{u}|^2 \, \mathrm{d}x \mathrm{d}t
		\le \int_{Q_2} h(E(|\dot{u}|)|\dot{u}|^2) \, \mathrm{d}x \mathrm{d}t.
	\end{align*}
	Now, Jensen's inequality yields
	\begin{align*}
		\int_{Q_2} h(E(|\dot{u}|)|\dot{u}|^2)\ \mathrm{d}x \mathrm{d}t \le |Q| h\left(
		\int_{Q_2} \frac{E(|\dot{u}|)|\dot{u}|^2}{|Q|} \, \mathrm{d}x \mathrm{d}t\right),
	\end{align*}
	and, recalling the definition of function $r$ in Equation (\ref{NLDefinOfR}), we get
	\begin{align}
		\notag
		\int_{Q_2} \big(E^2(|\dot{u}|)|\dot{u}|^2 + |\dot{u}|^2 \big) \mathrm{d}x \mathrm{d}t
		&\le |Q| r\left(\int_{Q_2} E(|\dot{u}|)|\dot{u}|^2 \, \mathrm{d}x \mathrm{d}t\right) \\
		\label{NLIneqForQ2} 
		&\le |Q| r\left(\int_{Q} E(|\dot{u}|)|\dot{u}|^2 \, \mathrm{d}x \mathrm{d}t\right).
	\end{align}
	Combining Equations (\ref{NLIneqForQ1}) and (\ref{NLIneqForQ2}), we obtain
	\begin{align*}
		\int_0^T \int_{\Omega} \big(E^2(|\dot{u}|)|\dot{u}|^2 + |\dot{u}|^2\big) \mathrm{d}x \mathrm{d}t \le  
		\left(c_1 \mathrm{I} + |Q|r\right) \left(\int_{Q} E(|\dot{u}|)|\dot{u}|^2 \,
		\mathrm{d}x \mathrm{d}t\right),
	\end{align*}
	where $\mathrm{I}$ is the identity function, i.e.,
	\begin{align*}
		\left(c_1 \mathrm{I} + |Q|r\right)(s) = c_1 s  + |Q| r(s) \text{ for } s \ge 0.
	\end{align*}

	Further, exploiting Lemma \ref{NLLemaEstimationForEnergy} and the monotonicity of $r$, we get
	\begin{align*}
		\nlE(T) &\le C \left(\left( c_1 \mathrm{I} + |Q|r\right) \left(\int_{Q} E(|\dot{u}|)|\dot{u}|^2 \, \mathrm{d}x \mathrm{d}t\right)
		+ \int_Q |q|^2 \ \mathrm{d}x \mathrm{d}t\right)
		\\
		&\le C \left(c_2 \int_{Q} \big(E(|\dot{u}|)|\dot{u}|^2 + |q|^2 \big) \mathrm{d}x \mathrm{d}t +
		|Q|r\left(\int_{Q} E(|\dot{u}|)|\dot{u}|^2 \, \mathrm{d}x \mathrm{d}t\right)\right)
		\\
		&\le C \left(c_2 \int_{Q} \big(E(|\dot{u}|)|\dot{u}|^2\ + |q|^2\big) \mathrm{d}x \mathrm{d}t + 
		|Q|r\left(\int_{Q} \big(E(|\dot{u}|)|\dot{u}|^2 + |q|^2\big) \mathrm{d}x \mathrm{d}t\right)\right)
		\\
		&= C \left((c_2 \mathrm{I} + |Q|r) \left(\int_{Q} \big(E(|\dot{u}|)|\dot{u}|^2  + |q|^2\big) \mathrm{d}x \mathrm{d}t\right)\right),
	\end{align*}
	where $c_2 = \max\{c_1, 1\}$. 
	Thus, we obtain
	\begin{align*}
		\frac{1}{C|Q|} \nlE(T) \le \left(\frac{c_2}{|Q|} \mathrm{I} + r\right)
		\left(\int_{Q} \big(E(|\dot{u}|)|\dot{u}|^2 + |q|^2\big) \mathrm{d}x \mathrm{d}t\right),
	\end{align*}
	or, equivalently,
	\begin{align*} 
		M \nlE(T) \le (c \mathrm{I} + r)(\nlE(0) - \nlE(T)).
	\end{align*}
	Now, recalling Equation (\ref{NLDefinOfP}), we get
	\begin{align*}
		p(\nlE(T)) \le \nlE(0) - \nlE(T),
	\end{align*}
	or, equivalently,
	\begin{align*}
		p(\nlE(T)) + \nlE(T) \le \nlE(0).
	\end{align*}
	It can easily be seen that, replacing $0, T$ with $mT, (m+1)T$ above, we obtain
	\begin{align*}
		p\big(\nlE((m + 1)T)) + \nlE((m + 1)T\big) \le \nlE(mT),
	\end{align*}
	Note that the constants $C, c, c_1, c_2, |Q|$ and the functions $h,r,p,q$ remain the same. 
	As before, applying Lemma \ref{LemmaLasiecka} with $s_m=\nlE(mT)$, 
	we conclude
	\begin{equation}
		E(mT) \le S(m) \text{ for } m \in \mathbb{N} \text{ and }
		\lim_{t\to\infty} S(t) = 0. \notag
	\end{equation}
	It can further be easily shown that $q(s) \ge 0$ for $s\ge 0$ and, thus,
	$S(\cdot)$ is monotonically decreasing. 
	Therefore, for $t = mT + \tau$ with $\tau \in [0, T]$,  we have
	\begin{align*}
		E(t) \le E(mT) \le S(m) = S\left(\frac{t-\tau}{T}\right) \le S\left(\frac{t}{T} - 1\right),
	\end{align*}
	which finishes the proof.
\end{proof}

\section*{Funding}
	This work has been partially funded by the Young Scholar Fund at the University of Konstanz, Germany.
	MP has been supported by the ERC-CZ Project LL1202 `MOdelling REvisited + MOdel REduction' at Charles University in Prague, Czech Republic and
	the Deutsche Forschungsgemeinschaft (DFG) through CRC 1173 at Karlsruhe Institute of Technology, Germany.

\appendix
\section{Relation between observability and stabilization}
\label{SECTION_APPENDIX}
Consider a pair of evolution equations
\begin{align}
	\label{APPEq1} 
	\dot{U}(t) + \oA U(t) + \oB U(t) &= 0 \text{ for } t \ge 0, & U(0) &= U_0,
	\\
	\label{APPEq2} 
	\dot{\phi}(t) + \oA \phi(t) &= 0 \text{ for } t \ge 0, & \phi(0) &= \phi_0,
\end{align}
taking place on a Hilbert space $\spH$ endowed with inner product $\langle\cdot, \cdot\rangle$. 
By $E_U(t) := \frac12 \langle U, U\rangle_\spH$ and $E_\phi(t) = \frac12 \langle \phi, \phi\rangle_\spH$ 
we denote the corresponding natural energies. 
Further, let $\oB \colon \spH \to \spH$ be a bounded linear operator. 
We are interested in establishing a relation between the following two properties: 
\begin{description}
        \item[$i)$]
        There exist numbers $T > 0$ and $C > 0$ such that every mild solution of Equation (\ref{APPEq2}) satisfies
        \begin{equation}
		\label{APPObserv}
		E_{\phi}(0) \le C \int_0^T \langle \oB \phi, \phi\rangle \mathrm{d}t.
        \end{equation}
        
        \item[$ii)$]
        There exist numbers $M > 0$ and $\lambda > 0$ such that every mild solution of Equation (\ref{APPEq1}) satisfies
        \begin{equation}
		\label{APPStabiliz}
		E_{U}(0) \le M e^{-\lambda t} E_U (0) \text{ for } t \geq 0.
        \end{equation}
\end{description}
Condition $i)$ is referred to the (exact) observability inequality at time $T$ for the operator pair $(\oA,\oB)$, 
whereas $ii)$ represents the exponential stability of the semigroup generated by $\oA + \oB$.
Using the approach proposed by Haraux \cite{Haraux1} and further developed Tebou \cite{Tebou1} 
for a class of abstract wave-equation-like problems, 
we prove the following theorem.
\begin{theorem}
	\label{APPTheorem}
	Assume $\oA$ is a infinitesimal generator of a $C_0$-semigroup of contractions on $\spH$ satisfying
	\begin{equation}
		\label{APPDissOfA}
		\langle\oA U, U\rangle = 0 \text{ for any } U \in D(\oA),
	\end{equation}
	and $\oB$ is a self-adjoint, positive semidefinite operator.
	Then there exits a number $T_{0} > 0$ such that
	the condition (\ref{APPStabiliz}) implies (\ref{APPObserv}) for any $T \geq T_{0}$.
\end{theorem}

\begin{proof}
	Let $\phi$ be an arbitrary strong solution of (\ref{APPEq2}). 
	By $U$ we denote a strong solution of the problem (\ref{APPEq1}) with the initial condition
	\begin{equation}
		U(0) = \phi(0). \notag
	\end{equation}
	Multiplying Equation (\ref{APPEq1}) by $U$ in $\spH$ and using Equation (\ref{APPDissOfA}), we get
	\begin{align*}
		\langle \dot{U}, U\rangle + \langle\oA U, U\rangle &= -\langle \oB U, U\rangle,
		\\
		\p E_{U}(t) &= -\langle\oB U, U\rangle,
		\\
		E_U(0) - E_U(T)  &= \int_0^T \langle\oB U, U\rangle \mathrm{d}t.
	\end{align*}
	It follows from Equation (\ref{APPStabiliz}) that there exists $T_{0} > 0$
	such that $E_U(T) \le \frac{1}{2} E_U(0)$ for $T \geq T_{0}$. Hence,
	\begin{equation}
		\int_0^T \langle\oB U, U\rangle \mathrm{d}t \ge \frac{1}{2} E_U(0). \notag
	\end{equation}
	Letting $\psi = \phi - U$, we observe $\psi$ is a strong solution of
	\begin{equation}
		\dot{\psi}(t) + \oA \psi(t) + \oB \psi(t) = \oB \phi(t) \text{ for } t > 0, \quad \psi(0) = 0. \notag
	\end{equation}
	Denoting by $E_{\psi}(t)$ the associated energy, we get
	\begin{align*}
		\langle\dot{\psi}, \psi\rangle + \langle\oA \psi, \psi\rangle 
		+ \langle\oB \psi, \psi\rangle &= \langle\oB \psi, \phi\rangle,
		\\
		\langle\dot{\psi}, \psi\rangle + 
		\langle\oB \psi, \psi\rangle &= \langle\oB \psi, \phi\rangle.
	\end{align*}
	Using Cauchy \& Schwarz' inequality for bilinear forms  (cf, e.g., \cite[p. 358]{Kadec1}), 
	we get the following cascade of inequalities
	\begin{align*}
		\langle\dot{\psi}, \psi\rangle + \langle\oB \psi, \psi\rangle 
		= \langle\oB \psi, \phi\rangle &\le \frac{1}{2}
		\left(\langle\oB \psi, \psi\rangle + 
		\langle\oB \phi, \phi\rangle\right),
		\\
		\langle\dot{\psi}, \psi\rangle + \frac{1}{2} \langle\oB \psi, \psi\rangle  
		&\le \frac{1}{2} \langle\oB \phi, \phi\rangle,
		\\
		\int_0^T \langle\dot{\psi}, \psi\rangle \mathrm{d}t + 
		\frac{1}{2} \int_0^T \langle\oB \psi, \psi\rangle \mathrm{d}t
		&\le \frac{1}{2} \int_0^T \langle\oB \phi, \phi\rangle \mathrm{d}t,
		\\
		E_{\psi}(T) - E_{\psi}(0) + \frac{1}{2} \int_0^T \langle\oB \psi, \psi\rangle \mathrm{d}t
		&\le \frac{1}{2} \int_0^T \langle\oB \phi, \phi\rangle \mathrm{d}t.
	\end{align*}
	Now, using the fact $E_{\psi}(0) = 0$ and $E_{\psi}(t)\ge 0$, we get
	\begin{align}
		\label{APPIntIneq}
		\int_0^T \langle \oB \psi, \psi\rangle \mathrm{d}t \le 
		\int_0^T \langle \oB \phi, \phi\rangle \mathrm{d}t.
	\end{align}
	Tacking into account the fact $E_U(0) = E_{\phi}(0)$, 
	Equation (\ref{APPIntIneq}) and Cauchy \& Schwarz' inequality for bilinear forms,
	we obtain
	\begin{align*}
		\frac{1}{2} E_\phi(0) &= \frac{1}{2} E_{U}(0) \le \int_0^T \langle\oB U, U\rangle \mathrm{d}t =
		\int_0^T \big\langle\oB (\phi-\psi), (\phi - \psi)\big\rangle \mathrm{d}t
		\\
		&= \int_0^T \langle\oB \phi, \phi\rangle \mathrm{d}t + 
		\int_0^T \langle\oB \psi, \psi\rangle \mathrm{d}t -
		2\int_0^T \langle\oB \phi, \psi\rangle \mathrm{d}t
		\\
		&\le 2\left(\int_0^T \langle\oB \phi, \phi\rangle \mathrm{d}t 
		+ \int_0^T \langle\oB \psi, \psi\rangle \mathrm{d}t\right) 
		\le 8 \int_0^T \langle\oB \phi, \phi\rangle \mathrm{d}t,
	\end{align*}
	which finishes the proof.
\end{proof}

\end{document}